\documentclass[sn-basic,iicol]{sn-jnl}
\usepackage{graphicx}%
\usepackage{multirow}%
\usepackage{amsmath,amssymb,amsfonts}%
\usepackage{amsthm}%
\usepackage{mathrsfs}%
\usepackage[title]{appendix}%
\usepackage{xcolor}%
\usepackage{textcomp}%
\usepackage{manyfoot}%
\usepackage{booktabs}%
\usepackage{algorithm}%
\usepackage{algorithmicx}%
\usepackage{algpseudocode}%
\usepackage{listings}
\usepackage{mathtools}
\usepackage[font=scriptsize]{subcaption}

\theoremstyle{thmstyleone}%
\theoremstyle{thmstyletwo}%

\theoremstyle{thmstylethree}%

\newcommand{\Ebb}{\ensuremath{\mathbb{E} }}
\newcommand{\Vbb}{\ensuremath{\mathbb{V} }}

\newcommand{\ignore}[1]{}

\begin{document}

\title[Constrained Multi-fidelity Bayesian Optimization]{A Constrained Multi-Fidelity Bayesian Optimization Method}

\author*[1]{\fnm{Jingyi} \sur{Wang}}\email{wang125@llnl.gov}

\author[1]{\fnm{Nai-Yuan} \sur{Chiang}}\email{chiang7@llnl.gov}
\author[1]{\fnm{Tucker} \sur{Hartland}}\email{hartland1@llnl.gov}

\author[1]{\fnm{J. Luc} \sur{Peterson}}\email{peterson76@llnl.gov}
\author[1]{\fnm{Jerome} \sur{Solberg}}\email{solberg2@llnl.gov}
\author[1]{\fnm{Cosmin G.} \sur{Petra}}\email{petra1@llnl.gov}

\affil[1]{\orgname{Lawrence Livermore National Laboratory}, \orgaddress{\city{Livermore}, \postcode{94550}, \state{CA}, \country{USA}}}

\abstract{
Recently, multi-fidelity Bayesian optimization (MFBO) has been successfully applied to many engineering design optimization problems, where the cost of high-fidelity simulations and experiments can be prohibitive. 
However, challenges remain for constrained optimization problems using the MFBO framework, particularly in efficiently identifying the feasible region defined by the  constraints. 
In this paper, we propose a constrained multi-fidelity Bayesian optimization (CMFBO) method with novel acquisition functions. 
Specifically, we design efficient acquisition functions that 1) have analytically closed-form expressions; 2) are straightforward to implement; and 3) do not require feasible initial samples --- an important feature often missing in commonly used acquisition functions such as expected constrained improvement (ECI). 
We demonstrate the effectiveness of our algorithms on synthetic test problems using different combinations of acquisition functions. 
Then, we apply the proposed method to a data-driven inertial confinement fusion (ICF) design problem, and a high-current joint design problem using finite element simulations with computational contact mechanics. }


\keywords{Bayesian optimization, constrained optimization, multi-fidelity, inertial confinement fusion, computational mechanics}



\maketitle

\newcommand{\norm}[1]{\left\lVert {#1} \right\rVert}
\newcommand{\Rbb}{\ensuremath{\mathbb{R} }}


\section{Introduction}\label{se:intro}
Bayesian optimization (BO) is a powerful, gradient-free, black-box optimization method.
It typically employs Gaussian process (GP) surrogate models to approximate the objective function using sample data, often obtained from physical experiments or numerical simulations. 
BO has been applied successfully to many areas, including structural design~\citep{mathern2021,ragueneau2024constrained,yoon2025application}, process optimization in additive manufacturing~\citep{wang2021coupled,wang2023optimization}, inertial confinement fusion (ICF) design~\citep{wang2023multifidelity}, hyper-parameter tuning in machine learning~\citep{wu2019hyperparameter,letham2019bayesian}. Recent advances have extended BO algorithms to more complex scenarios, including those with multi-fidelity data~\citep{zuluaga2013} and constraints~\citep{bernardo2011}.

Constrained Bayesian optimization (CBO) extends BO to constrained problems by modeling black-box constraint functions jointly with the objective function. The optimization problem it aims to solve can be represented as  
\begin{equation} \label{eqn:opt-prob}
 \centering
  \begin{aligned}
	  &\underset{\substack{x}\in C}{\text{minimize}} 
	  & & f(x), \\
          &\text{subject to}
	  & & c(x) \geq  0,
  \end{aligned}
\end{equation}
where $C\subset \Rbb^n$ defines bound constraints, $f:\Rbb^n\to \Rbb$ is the objective, and $c:\Rbb^n\to\Rbb^m$ represents the inequality constraints.
Various CBO methods have been proposed in the past decade, with the focus on the development of different acquisition functions~\citep{frazier2018}. 

One of the most popular CBO methods is the expected constrained improvement (ECI)~\citep{gelbart2014bayesian,gardner2014,letham2019constrained,eriksson2021scalable}. It builds on top of expected improvement (EI), which only requires evaluations of the cumulative distribution function (CDF) and probability density function (PDF) of the standard normal distribution.
Upper confidence bound (UCB) and its extension are another class of important acquisition functions that can be applied to constrained problems via penalty methods~\citep{lu2022no} or explicit constraint enforcement~\citep{zhou2022kernelized}.
Other notable CBO methods include the predictive entropy search with constraints (PESC)~\citep{hernandez2015predictive,takeno2022sequential}, augmented Lagrangian (AL) methods~\cite{gramacy2016modeling}, Slack-AL~\citep{picheny2016bayesian}, ADMMBO~\citep{ariafar2019admmbo}, active learning feasible region search~\citep{zhang2023constrained,jiao2019complete}, etc.

On the other hand, multi-fidelity Bayesian optimization (MFBO) has also garnered significant interest in recent years, particularly for complex engineering systems with expensive high-fidelity simulations and experiments. By leveraging strongly correlated low-fidelity models, MFBO can potentially achieve accurate solutions at a reduced cost. Notable applications include optimal designs for drone~\citep{charayron2021multi}, inertial confinement fusion (ICF)~\citep{wang2023multifidelity}, and airfoil shape design~\citep{meliani2019multi}. 
Multi-fidelity surrogate models, particularly multi-fidelity cokriging~\citep{raissi2016deep,raissi2017machine,perdikaris2015multi,ghoreishi2019multi,xiao2018extended,legratiet2013} have been widely studied and adopted in practice. 
Among MFBO methods, EI is one of the most popular acquisition functions~\citep{zhang2018variable,sarkar2019,shu2021multi}. 
A more thorough review of MFBO can be found in~\cite{do2023multi}, where the authors conclude that constrained MFBO is a topic for future research. 

Many engineering optimization problems are subjected to black-box constraints that can be expensive to evaluate. This has led to recent efforts to develop constrained multi-fidelity Bayesian optimization (CMFBO) methods.  
Despite considerable progress, developing practical CMFBO algorithms remains challenging due to the limited availability of acquisition functions.
Most existing CMFBO approaches rely on ECI and the probability of feasibility (PoF) used in ECI. 
However, ECI and PoF require feasible initial samples, a condition not guaranteed in practice~\citep{gardner2014}. 
Furthermore, in the MFBO framework, ECI can suffer from poor exploration of the boundary of the feasible region~\citep{zhou2024multi}.
Finally, it is reported in literature that the performance of ECI still has room for improvement \citep{hernandez2015predictive,ariafar2019admmbo,picheny2016bayesian}.
On the other hand, high-performing CBO methods such as PESC~\citep{hernandez2015predictive} have yet to be widely applied to multi-fidelity engineering design optimization problems, partially because of the lack of closed-form expressions, increased algorithm complexity, and higher implementation cost.

As mentioned above, the widely adopted ECI acquisition function is a natural extension of EI to the constrained setting, where the EI of the objective function is multiplied by the PoF function of the constraints.
CMFBO methods that use ECI as acquisition functions are presented in~\cite{khatouri2020constrained,sarkar2019,tran2020smf,shu2021multi,zhou2024multi}. 
In~\cite{sarkar2019}, the authors apply EI for building high-fidelity sampling, and mutual information (MI) for low-fidelity sampling for unconstrained problems. They discussed using PoF and ECI to tackle constraints in their MFBO framework. In~\cite{shu2021multi}, the authors proposed an expected further improvement function for MFBO with a hierarchical kriging model, and extended it to the constrained setting using the ECI approach. In~\cite{tran2020smf}, a variation of the PoF is computed from GP surrogate models, referred to as a binary classifier surrogate model. In~\cite{zhou2024multi}, the authors identified similar challenges with existing CMFBO methods and proposed an approach that does not solely use ECI. Their solution is a three-stage optimization method, in which the first stage focuses on finding a feasible initial sample and the third stage uses a weighted lower confidence bound (LCB) function to select an additional batch of samples. While this approach addresses the issue of infeasible initial samples, the algorithm focuses on multi-objective problems and can be less efficient due to its exclusive focus on feasibility in the first stage.

In this paper, we develop a cokriging-based CMFBO method  through multiple novel acquisition functions with closed-form expressions that are effective, do not require feasible initial samples, and are simple to implement. 
Our new acquisition functions are derived from the classic penalty methods in nonlinear optimization~\citep{Nocedal_book}. In particular, we incorporate constraint violations directly into the objective, offering both simplicity and ease of implementation. To the best of the authors' knowledge, this work is among the first to present a CMFBO method that does not rely on the standard ECI, and thus addresses the challenges associated with applying ECI.
We validate our method through synthetic benchmark problems and two engineering design problems based on real-world applications: a data-driven ICF design problem and a high-current joint design problem modeled by computational contact mechanics. 
In addition, we study the impact of the number of low-fidelity samples and the choice of aquisition functions on the convergence performance of the method. 
The main contributions of this paper can be summarized as follows.
\begin{itemize}
  \item We propose the expected merit improvement (EMI), a closed-form acquisition function, and an update rule for its algorithmic parameter. Further, we develop the additive expected constrained
improvement (AECI) that combines EMI with ECI.
  Both EMI and AECI use additive structures between the objective and the constraints, and do not need feasible initial samples. 
  \item We design the constrained upper confidence bound (CUCB) acquisition function by incorporating the penalty constraint violation function. CUCB is easy to implement, has closed-form expressions, and do not require feasible initial samples.
  \item We present our multi-fidelity strategy that is generalizable and applicable to different acquisition functions. We conduct tests on the effect of using different  acquisition functions for low-fidelity sampling under the proposed method. 
  \item We apply the proposed method to two engineering applications.
  For the ICF design problem, we use simulation data generated from the multi-physics simulation code HYDRA~\citep{metal2001}. 
  For the current-joint design problem, we model it using finite element method and mortar contact formulation, and run the simulation as the optimization progresses.  
  To the best of our knowledge, this paper is among the first to apply BO methods to the optimization of objects in contact, as well as constrained ICF design problem. 
  Thus, we validate the use of BO methods for both applications.    
\end{itemize}
The remainder of this paper is organized as follows.  Section~\ref{se:cbo} provides background on GP and CBO methods, while introducing our proposed acquisition functions. Section~\ref{se:multi-bayesian} presents the multi-fidelity strategy and the proposed CMFBO algorithm. Section~\ref{sec:exp} presents numerical experiments on synthetic and two real-world design optimization problems. Section~\ref{sec:con} gives a conclusion of this paper.

\section{Constrained Bayesian optimization}\label{se:cbo}
\par\noindent
In this section, we briefly introduce two main components  of a CBO algorithm: the GP surrogate models for both the objective function $f$ and the constraint functions $c_i,i=1,\dots,m$, and the acquisition function that is used to select the next samples~\citep{frazier2018,bayesianoptreview2016}. 
Then, we present our novel acquisition functions.

\subsection{Gaussian process}\label{se:cbo-bayesian}
The GP model defines a joint multivariate Gaussian distribution over function evaluations at different inputs $x$. 
A GP is specified by a mean function $m(x):\Rbb^n\to\Rbb$ and a covariance function (kernel) $k(x,x'):\Rbb^n\times\Rbb^n\to\Rbb$.
In practice, the mean function $m(\cdot)$ is often chosen to be the zero constant function $m(x)\equiv 0$, which is also adopted in 
 this paper. The zero-mean GP model is denoted as $GP(0,k(x,x'))$.
 Given $t$ input samples $x_1, \dots, x_t$, the corresponding objective observations can be written as $y_i = f(x_i) + \epsilon_i$, $i=1,\dots,t$,  where the noise is assumed to follow a zero-mean Gaussian distribution, \textit{i.e.}, $\epsilon_t \sim \mathcal{N}(0,\sigma^2)$.
 Denote the observations as $y_{1:t}=[y_1,\dots,y_t]^T$ at $x_{1:t}=[x_1,\dots,x_t]^T$.
 The posterior mean and variance of the GP can be inferred to be 
\begin{equation} \label{eqn:GP-post}
 \centering
  \begin{aligned}
  \mu_t(x) =&~ k(x,x_{1:t}) [K(x_{1:t},x_{1:t})+ \sigma^2 I]^{-1} y_{1:t} ,\\
  \sigma^2_t(x) =&~ k(x,x)-k(x_{1:t},x)^T[K(x_{1:t},x_{1:t})\\
                 &+\sigma^2 I]^{-1} k(x_{1:t},x).
\end{aligned}
\end{equation}
Here, $K(x_{1:t},x_{1:t})$ is the covariance matrix and $k(x_{1:t},x)=[k(x_1,x),\dots,k(x_t,x)]^T$. 
The full posterior covariance kernel can be derived as 
$
k_t(x, x') = k(x, x') - k(x_{1:t},x)^T [K(x_{1:t}, x_{1:t}) + \sigma^2 I]^{-1} k(x_{1:t}, x')
$. 
Both the objective and constraint functions are approximated by their respective GP surrogate models.  
We can similarly obtain the posterior mean and variance of the GP for each constraint function.
We use the superscripts $^f$ and $^{c_j}$, $j=1,\dots,m$ to differentiate the GP models and posterior predictions for the objective and constraint functions.
That is, $\mu^{c_j}_t(x)$ and $(\sigma^{c_j}_t(x))^2$ denote the GP posterior mean and variance for constraint function $c_j(x)$, while  $\mu^{f}_t(x)$ and $(\sigma^{f}_t(x))^2$ denote the GP posterior mean and variance for the objective.
For simplicity of presentation, we do not add additional notations to distinguish between the kernels for the objective and constraints. However, we emphasize that they can be different functions.

The choice of kernel $k(\cdot,\cdot)$ has significant impact on the performance of the BO algorithms. The squared exponential (SE) 
 and Matérn kernels are among the most popular kernels. The SE kernel is defined as follows
\begin{equation} \label{def:se}
  \centering
  \begin{aligned}
    k_{SE}(x,x';\theta) =& \exp\left(-\frac{r^2}{2\theta^2}\right), \\
  \end{aligned}
\end{equation}
where $\theta>0$ is the length hyper-parameter and $r=\norm{x-x'}_2$ for $x,x'\in C$. 
 The hyper-parameters of a GP can be typically estimated by maximizing the log-marginal-likelihood of the training data with an optimization method such as the L-BFGS algorithm~\citep{Nocedal_book}. 

\subsection{Expected improvement and upper confidence bound}
Acquisition functions determine the next sample points at each BO iteration. 
Two of the most popular acquisition functions for unconstrained BO are EI and UCB, which form the basis for our novel acquisition functions in constrained setting. 

EI is one of the most successful and widely used acquisition functions~\citep{brochu2010}. It takes the conditional expectation  $\Ebb_{t}$ over $t$ samples of the improvement function defined as
\begin{equation} \label{eqn:improvement}
    I_t(x) = \max\{ y^+_{t} -f(x),0  \}, 
\end{equation}
where 
     $y^+_{t} = \min\{y_1, \dots, y_t\}$
is the best observed objective value so far.
From~\eqref{eqn:improvement}, EI has the following closed-form expression:
\begin{equation} \label{eqn:EI-1}
EI_t(x) = z_t(x) \sigma_t(x) \Phi(z_t(x)) + \sigma_t(x) \phi(z_t(x)),
\end{equation}
where $z_t(x) = \frac{y_t^+ - \mu_t(x)}{\sigma_t(x)}$, and $\phi(\cdot)$ and $\Phi(\cdot)$ are the standard normal PDF and CDF, respectively. 
The next sample point is chosen by maximizing EI, \textit{i.e.},
\begin{equation} \label{eqn:acquisition-1}
x_{t+1} = \arg\max_{x \in C} EI_t(x).
\end{equation}

The UCB acquisition function balances exploration and exploitation by directly using an algorithmic parameter, the posterior mean, and variance of the GP model~\citep{srinivas2009gaussian}. 
UCB and its extension such as mutual information (MI)~\citep{contal2014} have been reported to be effective for selecting low-fidelity samples due to its controllable exploration properties~\citep{contal2014,sarkar2019}. 
The mathematical formulation of UCB for maximization problems is 
\begin{equation} \label{eqn:GPucb}
UCB_t(x) = \mu_t(x) + \sqrt{\beta_{t+1}} \sigma_t(x),
\end{equation}
where $\beta_{t+1} > 0$ is a parameter that balances exploration and exploitation. 
Similar to EI, the next sample point is chosen by maximizing the UCB function 
\begin{equation} \label{eqn:acquisition-ucb}
x_{t+1} = \arg\max_{x \in C} UCB_t(x).
\end{equation}
The choice of the parameter $\beta_t$ has been studied in the literature and often takes into consideration the theoretical cumulative regret properties of UCB~\citep{srinivas2009gaussian,chowdhury2017kernelized}. 

\subsection{Expected constrained improvement}\label{se:cbo-eci}
ECI~\citep{bernardo2011,gelbart2014bayesian,gardner2014} extends EI to constrained problems by incorporating the PoF function constructed from the constraint.
The constrained improvement function is defined as  
\begin{equation} \label{eqn:c-improvement}
 \centering
  \begin{aligned}
    I_{t}^C(x) = \Delta_{t}(x)\max\{ y^+_{t} -f(x),0  \} = \Delta_{t}(x) I_t(x), 
  \end{aligned}
\end{equation}
where $\Delta_t(x)\in\{0,1\}$ is a feasibility indicator function, taking the value $1$ if $x$ is feasible and $0$ otherwise. 
The best feasible objective from $t$ samples is
\begin{equation} \label{eqn:c-improvement-2}
    y_t^+ = \min\{ y_i \mid 1 \leq i \leq t, \; c_i(x) \geq 0 \}.
\end{equation}
For $y_t^+$ to be well-defined, at least one feasible (initial) point must be observed.
If the objective and constraint functions are conditionally independent, then the ECI function has the closed-form expression
\begin{equation} \label{eqn:cEI-1}
       EI^C_{t}(x) = PF_t(x) EI_t(x),
\end{equation}
where $PF_t(\cdot)$ is the probability that $x$ is feasible, \textit{i.e.}, the PoF function, and it follows a univariate Gaussian CDF when there is a single constraint or $m=1$. 

For problems with multiple constraints, they are assumed to be conditionally independent of each other in ECI. Consequently, the PoF is given by 
\begin{equation} \label{eqn:cEI-2}
     PF_t(x) = \prod_{j=1}^{m} P(c_j(x)\geq 0), 
\end{equation}
where $c_j$ denotes the $j$th constraint and $P$ is the probability operator. 
The next sample is selected by maximizing the ECI acquisition function: 
\begin{equation} \label{eqn:acquisition-ECI}
 \centering
  \begin{aligned}
      x_{t+1} = \underset{\substack{x\in C}}{\text{argmax}}  \,EI^C_{t} (x).
  \end{aligned}
\end{equation}
%

\subsection{ Expected Merit Improvement} \label{se:merit}
In this section, we propose an acquisition function using the penalty-based merit functions commonly used in nonlinear optimization~\citep{Nocedal_book}, which is designed to overcome the challenges associated with ECI, as outlined in Section~\ref{se:intro}. The exact penalty merit function for problem~\eqref{eqn:opt-prob} is defined as
\begin{equation} \label{eqn:merit-function}
 \centering
  \begin{aligned}
    \varphi_t(x) :=& f(x) + \alpha_t \sum_{j=1}^m c^+_j (x),
  \end{aligned}
\end{equation}
where $\alpha_t > 0$ is a penalty parameter, and $c_j^+(x):= \max\{-c_j(x), 0\}$ defines the violation of the $j$-th constraint. 
During the optimization process, the penalty parameter often needs to be updated to ensure that the constraint is sufficiently enforced. Thus, at $t$ samples, the corresponding $\alpha_t$ is dependent on $t$. 
We use a different symbol $\psi$ to denote the observed merit function value, similar to $y_i$ being an observed objective function value.
At a sample $x_i, i=1,\dots,t$, the observed merit function value is dependent on $t$ in the following form
\begin{equation}
  \psi_i(t) = y_i + \alpha_t \sum_{j=1}^m c^+_j(x_i),
\end{equation}
where $y_i$ is the observed objective value. The observed merit function values is written as $\psi_{1:t}(t) = [\psi_1(t), \dots, \psi_t(t)]$.

The penalty parameter $\alpha_t$ can be set to a sufficiently large constant to ensure feasibility, or updated adaptively following some standard techniques from nonlinear optimization~\citep{Nocedal_book}. We propose a practical update rule for $\alpha_t$ at the end of this section in Algorithm~\ref{alg:penalty}.

Using the penalty merit function~\eqref{eqn:merit-function}, we can extend the EI approach to the constrained setting. 
To do so, one needs to define an improvement function of $\varphi_t$ and take its conditional expectation. However, since $\varphi_t$ itself does not follow multivariate Gaussian distribution anymore (due to $c_j^+(x)$), we do not have closed-form expressions for the expectation of its improvement function, which now includes nested max functions. 
Thus, this approach could add considerable computational cost and increase implementation complexity.

To overcome this challenge, we propose the novel   merit improvement function as follows 
\begin{equation} \label{eqn:m-improvement-1}
  \begin{aligned}
    I_{t}^M(x) =& \max\{ y_{t^+}-f(x),0\} \\
              &+ \alpha_t \sum_{j=1}^m [c_j(x_{t^+}) -  c^+_j(x)],\\ 
  \end{aligned}
\end{equation}
where the subscript $t^+\in \{1,2,\dots,t\}$ is an index and used to distinguish~\eqref{eqn:m-improvement-1} from the improvement function $I_t(x)$. 
The index $t^+$ is defined as  
\begin{equation} \label{eqn:m-max-index}
  \begin{aligned}
     t^+=  \underset{\substack{i=1,\dots,t}}{\text{argmax}} \  \psi_t(x_i).
  \end{aligned}
\end{equation}
  Using notations from Section~\ref{se:cbo-bayesian}, the expected value of~\eqref{eqn:m-improvement-1}, referred to as the \emph{expected merit improvement} (EMI), has the following closed form
\begin{equation} \label{eqn:m-ei-1}
  \begin{aligned}
    EI^M_{t}(x) &= z^f_t(x)\, \sigma^f_t(x)\, \Phi(z^f_t(x)) 
               + \sigma^f_t(x)\, \phi(z^f_t(x)) \\
               &+ \alpha_t c_{t}^m 
               + \alpha_t \sum_{j=1}^m [
                    \mu_t^{c_j}(x)\, \Phi(z_t^{c_j}(x)) 
                    \\
                    &- \sigma_t^{c_j}(x)\, \phi(z_t^{c_j}(x)) 
               ],
  \end{aligned}
\end{equation}
where
\begin{equation} \label{eqn:m-ei-2}
  \begin{aligned}
    c_{t}^m &:= \sum_{j=1}^m c_{j}(x_{t^+}), \\
    z^f_t(x) &:= \frac{y_{t^+} - \mu^f_t(x)}{\sigma^f_t(x)}, \\
    z_t^{c_j}(x) &:= \frac{-\mu_t^{c_j}(x)}{\sigma_t^{c_j}(x)}, \quad j = 1,\dots,m.
  \end{aligned}
\end{equation}

Next, we  present our update rule for $\alpha_t$ in Algorithm~\ref{alg:penalty}. 
First, an initial value of $\alpha_0$ is chosen based on the knowledge of the problem. Then, at each iteration, Algorithm~\ref{alg:penalty} is called.  
We find the current best merit function value $\psi^+_t = \psi_t(x_{t^+})$ using~\eqref{eqn:m-max-index}.
If $\psi^+_t$ corresponds to a feasible sample point, $\alpha_t$ does not change. If $\psi^+_t$ is from an infeasible sample point, then $\alpha_t$ is increased via a fixed ratio $c_{\alpha}>1$.
Intuitively, a larger penalty puts more emphasis on the constraint, leading the algorithm eventually to a large enough $\alpha$ to effectively enforce the constraint.

\begin{algorithm}[t]
 \caption{Penalty update rule at iteration $t$}\label{alg:penalty}
  \begin{algorithmic}[1]
	  \State{Choose the increase ratio $c_{\alpha}>1$.} 
	  \State{Find the current maximum merit function value $\psi_t^+=\psi_t(x_{t^+})$.}
        \For{$j=1,\dots,m$ }
	    \If{$c_j(x_{t^+})<0$}
    	  \State{$\alpha_{t+1} = c_{\alpha}\alpha_t$. \;}
          \State{Break.\;}
	      \Else
            \State{$\alpha_{t+1} = \alpha_t$.\;}
          \EndIf
	\EndFor
  \end{algorithmic}
\end{algorithm}

\subsection{Additive expected constrained improvement}

The additive structure between the objective and constraint functions of EMI enables a natural way to combine it with ECI through an algorithmic weighting parameter $\beta\in[0,1]$. This hybrid acquisition function can be particularly beneficial for problems where ECI performs well and a sufficient number of feasible samples are available.
We define the following improvement function
\begin{equation} \label{eqn:u-improvement-1}
  \begin{aligned}
    I_{t}^{A}(x) := (1 - \beta)\, \Delta_t(x)\, I_t(x) + \beta\, I_{t}^M(x),
  \end{aligned}
\end{equation}
where the superscript $^A$ denotes ``additive''.
Its corresponding expected value, referred to as the additive expected constrained improvement (AECI) acquisition function, is given by
\begin{equation} \label{eqn:ceci-1}
 \centering
  \begin{aligned}
    EI_{t}^A  =& (1-\beta) EI^C_{t}(x) + \beta EI^M_{t}(x).\\ 
  \end{aligned}
\end{equation}
The algorithmic parameter $\beta \in [0,1]$ controls the balance between ECI and EMI. When $\beta = 0$, AECI~\eqref{eqn:ceci-1} reduces to ECI~\eqref{eqn:cEI-1}; when $\beta = 1$, it reduces to EMI~\eqref{eqn:m-ei-1}. Thus, both the classical ECI approach and our proposed EMI approach are unified under this single function. Since AECI inherits a closed-form expression from ECI and EMI, it continues to enable efficient and straightforward implementation.

In this paper, we use a simple rule based on feasibility to choose the parameter $\beta$. Specifically, let $n_f$ be the number of feasible samples observed so far, and let $N_f$ be a threshold predetermined by the user.  We then set
$$
  \beta = 
  \begin{cases}
    1, & \text{if } n_f < N_f, \\
    0, & \text{otherwise}.
  \end{cases}
$$
This rule prioritizes EMI during the early stages of the optimization, when few feasible points are available, and transitions to ECI once the feasible region has been sufficiently explored according to the practitioners' preference.  More sophisticated strategies for selecting $\beta$ can also be applied, such as employing a continuous value that balances ECI and EMI adaptively. A detailed investigation of such strategies is beyond the scope of this paper and a topic left for future work.

\subsection{Constrained upper confidence bound}\label{se:cucb}
In this section, we extend UCB to the constrained setting via the merit function~\eqref{eqn:merit-function}. First, we  apply the UCB formulation on $\varphi_t$, which leads to the form
  \begin{equation}\label{eqn:cucb-1}
    \begin{aligned}
	     & \Ebb_t [\varphi_t(x)] + \sqrt{\beta_{t+1}}\sqrt{\Vbb_t[\varphi_t(x)]}=\mu_t^f(x)\\
        &+\alpha_t\Ebb_{t}[c_{j}^+(x)] +\sqrt{\beta_{t+1}}\sqrt{\Vbb_t[\varphi_t(x)]} ,
    \end{aligned}
  \end{equation}
where $\Vbb_{t}$ denotes the conditional variance on $\varphi_t$. 
 Due to the non-linearity of the variance operator,~\eqref{eqn:cucb-1} remains complex to compute. We simplify it to the following form
 \begin{equation}\label{eqn:cucb-2}
    \begin{aligned}
	    \mu_t^f(x)+&\alpha_t\Ebb_{t}[c_{j}^+(x)] +\sqrt{\beta_{t+1}} \left(\sigma_t^f(x)\right.\\
        &\left.+\alpha_t \sum_{j=1}^m \sqrt{\Vbb_t[c_{j}^+(x)}]\right),
    \end{aligned}
  \end{equation}
  where we consider the variance of the objective and constraint components separately. 
The expectation and variance of $c_j^+$ have the following closed forms 
  \begin{equation}\label{eqn:cucb-3}
    \begin{aligned}
	   \Ebb_{t}[c_{j}^+(x)] =& -\mu^{c_j}_{t}(x)\Phi(z^{c_j}_{t}(x)) +\sigma^{c_j}_{t} \phi(z^{c_j}_{t}(x)),\\ 
           \Vbb_{t}[c_{j}^+(x)] =& \Ebb_{t}[(c_j^+(x))^2] - \Ebb^2_{t}[c_j^+(x)],\\
          \Ebb_{t}[(c^{+}_{j}(x))^2] =&   ((\mu^{c_j}_{t}(x))^2+(\sigma^{c_j}_{t}(x))^2)\Phi(z^{c_j}_{t}(x))\\
                      &- \mu^{c_j}_{t}(x) \sigma^{c_j}_{t}(x)\phi(z^{c_j}_{t}(x))],
    \end{aligned}
  \end{equation}
where  $z^{c_j}_{t}$ is from~\eqref{eqn:m-ei-2}.
While~\eqref{eqn:cucb-2} and~\eqref{eqn:cucb-3} return a closed-form expression, $\sqrt{\Vbb_{t}[c^{+}_{j}(x)]}$ as an uncertainty measure can be further simplified to avoid potential numerical issues arising from  complex expressions. For instance, the right-hand side of $\Vbb_{t}[c^{+}_{j}(x)]$ in~\eqref{eqn:cucb-3} can return small negative values when its true value is close to $0$ due to numerical errors.
As a result, our proposed CUCB acquisition function is 
 \begin{equation}\label{eqn:cucb-4}
    \begin{aligned}
	    &\mu_t^f(x)+\alpha_t\Ebb_{t}[c_{j}^+(x)] \\&+\sqrt{\beta_{t+1}} \left(\sigma_t^f(x)+\alpha_t\sum_{j=1}^m \sigma^{c_{j}}_t(x)\right),
    \end{aligned}
  \end{equation}
where $\Ebb_{t}[c_{j}^+(x)]$ is from~\eqref{eqn:cucb-3}.
We note that~\eqref{eqn:cucb-4} is derived from the UCB form for maximization problems. The minimization problem~\eqref{eqn:opt-prob} is equivalent to maximizing $-f$. Hence, we replace the posterior mean terms in~\eqref{eqn:cucb-4} with their negative counterparts when solving~\eqref{eqn:opt-prob}.
%
%

\section{Constrained multi-fidelity Bayesian optimization}\label{se:multi-bayesian}

Multi-fidelity Bayesian optimization  methods leverage simulations or experiments of varying fidelity to build surrogate models, offering efficiency gains when lower-fidelity simulations are significantly cheaper but still capture the key trends of higher-fidelity outputs. For constrained problems, multi-fidelity surrogate models are built not only for the objective but also for the constraints. In this section, we propose a cokriging-based CMFBO method, using the acquisition functions introduced in Sections~\ref{se:cbo}.

Let $L$ denote the number of fidelity levels, with the superscript $0$ indicating the lowest fidelity. We adopt the widely used additive multi-fidelity model~\citep{kennedy2001bayesian} in which each fidelity level is treated as a GP, and the difference between two adjacent fidelities is also modeled as a GP. Specifically, the surrogate model for the objective function at $i$th fidelity level is
\begin{equation}\label{eqn:multi-obj}
  \begin{aligned}
    f^0(x) &\sim GP(0, k(x,x')), \\
    f^i(x) &= \rho^{i-1} f^{i-1}(x) + \delta^{i-1}(x), \ 1 \leq i < L,
  \end{aligned}
\end{equation}
where $\rho^{i-1}$ is a scalar hyper-parameter, and $\delta^{i-1}(x)$ is a GP representing the gap between fidelities $i$ and $i-1$. Although this approach introduces additional surrogate models for each fidelity, it provides information that may be useful in post-processing.
Similarly, the constraint functions at different fidelities can be represented as 
\begin{equation}\label{eqn:multi-cons}
  \begin{aligned}
    c_j^0(x) &\sim GP(0, k(x,x')), \\
    c_j^i(x) &= \tau_j^{i-1} c_j^{i-1}(x) + \delta_j^{i-1}(x), \  1 \leq i < L,
  \end{aligned}
\end{equation}
where $c_j^i(x)$ denotes the $j$th constraint at fidelity level $i$; $\tau_j^{i-1}$ is a scalar; and $\delta_j^{i-1}(x)$ is again a standard GP model.

For clarity, we describe our CMFBO algorithm in the case of $L=2$, using superscripts $l$ and $h$ to denote the low- and high-fidelity levels, respectively. 
The predictive distributions for the two-fidelity model can be derived using the cokriging approach~\citep{legratiet2013}.
Each surrogate model is determined by its own posterior mean and standard deviation, scalar scaling parameters $(\rho^{l}, \tau^l)$, and kernel hyper-parameters, \textit{e.g.},  $(\theta^l, \theta^h)$ for the low- and high-fidelity objectives using the SE kernel~\eqref{def:se}. To improve efficiency, we estimate the posterior predictions and the hyper-parameters via Bayesian estimation in~\cite{legratiet2013} rather than a fully Bayesian inference procedure. 

The novel acquisition functions in Section~\ref{se:cbo} are used in the proposed CMFBO method, which incorporates flexible choices of acquisition functions at different fidelity levels.
While using different acquisition functions for multi-fidelity Bayesian optimization algorithms is not new~\citep{sarkar2019}, the increased options provided in Section~\ref{se:cbo} allow us to deploy more combinations of acquisition functions in the constrained setting to satisfy different application requirements and fidelity costs.
 For example, the CUCB acquisition function with a large value of $\beta_t$ can be applied to low-fidelity models to encourage exploration. 

The cokriging-based CMFBO algorithm for two-fidelity levels is outlined in Algorithm~\ref{alg:cmfbo}. At each iteration, we ensure $x_{1:t}^h \subseteq x_{1:t}^l$, meaning all the high-fidelity samples are also included in the low-fidelity sample set to maintain consistency across fidelities in the cokriging process. 
Moreover, each high-fidelity simulation may be accompanied by multiple low-fidelity simulations in the same iteration. A user-defined integer $N_\ell$ controls the ratio of additional low- to high-fidelity samples per iteration. This parameter can be tuned according to available computational resources. For example, when the simulation costs for each fidelity are known, $N_\ell$ can be chosen to synchronize the completion times of high- and low-fidelity simulations at each iteration~\citep{pmlr-v119-takeno20a,tran2020smf}. If $N_\ell = 0$, no additional low-fidelity samples are added, which implies that only one low-fidelity sample can be added in the cokriging setting.

\begin{algorithm}[t]
 \caption{Constrained multifidelity Bayesian optimization}\label{alg:cmfbo}
  \begin{algorithmic}[1]
	  \State{Choose initial sample points $x_{0}^{l}$, $x_{0}^h$ for low- and high-fidelity simulation.} 
	  \State{Build initial multi-fidelity surrogate models.}
    \For{$t=1,2,\dots$}
	  \State{Obtain the high-fidelity sample point $x_{t}^h$ by optimizing a high-fidelity acquisition function, e.g., AECI~\eqref{eqn:ceci-1}.}
	  \State{Run high-fidelity and low-fidelity experiments at $x_{t}^h$ to obtain objective and constraint observations.\;}
	  \State{Retrain GP models with the new samples.\;}
        \For{$j=0,1,\dots,N_\ell$}
	    \State{Evaluate the low-fidelity acquisition function, e.g., AECI. Find the new low-fidelity sample point $x_{t,j}^l$.}
    	  \State{Run low-fidelity experiments at $x_{t}^l$, obtaining objective and constraint observations. \;}
	      \State{Retrain GP models with the new samples.\;}
	\EndFor
	  \If{Stopping criteria satisfied} {Exit.}
          \EndIf 
  \EndFor
  \end{algorithmic}
\end{algorithm}

\newcommand\rhobold{{\ensuremath{\boldsymbol{\rho}}}}
\section{Numerical Experiments}\label{sec:exp}
In this section, we present synthetic test problems and two application-based constrained  optimization problems, and use the proposed CMFBO algorithms with EMI, AECI, CUCB, and ECI acquisition functions to solve them. 
Algorithm~\ref{alg:cmfbo} is implemented in Python using scikit-learn~\citep{scikit-learn,sklearn_api}, OpenMDAO~\citep{Gray2019a}, and PyTorch~\citep{paszke2019pytorch},  while the ICF and contact mechanics simulation codes are implemented in Python and MATLAB, respectively.  

Our numerical experiments showcase the following three key aspects of the algorithmic performance of the proposed CMFBO method.
First, we demonstrate the effectiveness of the proposed acquisition functions, where the algorithm succeeds in finding sufficiently accurate optimal solutions quickly in all examples. 

Second, we study the effect of the number of low-fidelity samples on the performance of the algorithm using  synthetic problems 1 to 8. Noticeably, we design the constraints in example 3 and 4  such that the high-fidelity and low-fidelity constraints have low correlation, while the rest of the synthetic examples have highly correlated constraints. 

Third, we compare the proposed CMFBO method to existing  widely used optimization methods in our application-based examples.
For the ICF design problem, we implement the commonly used ECI acquisition function. To obtain a fair comparison, we only change the acquisition function in CMFBO to ECI for both high and low fidelities. 
For the joint/wedge design problem, we compare the CMFBO method to both the ECI approach and a state-of-the-art gradient-based optimization algorithm. 

Unless otherwise specified, we choose a small increase ratio $c_{\alpha}=1.1$ for $\alpha_t$ when using EMI and AECI. Further, we set $N_f=2$ for AECI. The parameter $\beta_t$ is chosen to be constant $1$ for simplicity. The penalty parameter $\alpha_t$ is updated according to Algorithm~\ref{alg:penalty}.

\subsection{Synthetic examples}\label{se:ex-synthetic}
In this section, we present numerical experiments on synthetic constrained multi-fidelity optimization problems. 
The objective functions are chosen from unconstrained multi-fidelity optimization problems commonly found in the literature. Then, we incorporate distinct high-fidelity and low-fidelity constraints to construct constrained problems. 
We run Algorithm~\ref{alg:cmfbo} with Latin hypercube sampling for the initial samples of each problem. To obtain  statistically sound results, each problem is solved 100 times. 
Altogether, eight examples are considered, differentiated by the objective, constraint functions, and the acquisition functions used. 
Three levels of low-fidelity samples ($N_\ell = 0, 1, 2$) are tested in order to demonstrate the effect of low-fidelity sampling on optimization performance. Five initial samples are used unless otherwise specified.

\begin{figure*}[ht]
    \begin{subfigure}{0.48\textwidth}
        \centering
	\includegraphics[width=\linewidth]{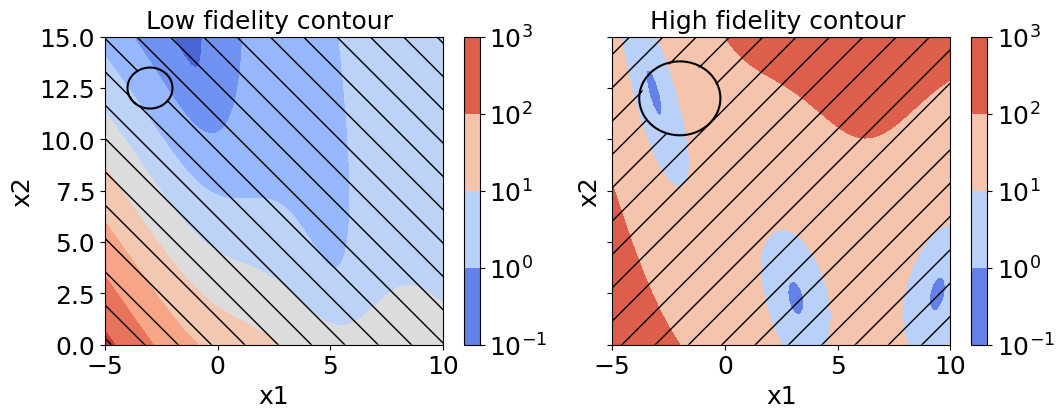}
    \caption{Low- and high-fidelity contour plots of example 1 and 2.}
  \label{fig:ex-Branin}
    \end{subfigure}
    \hfill
    \begin{subfigure}{0.48\textwidth}
        \centering
	\includegraphics[width=\linewidth]{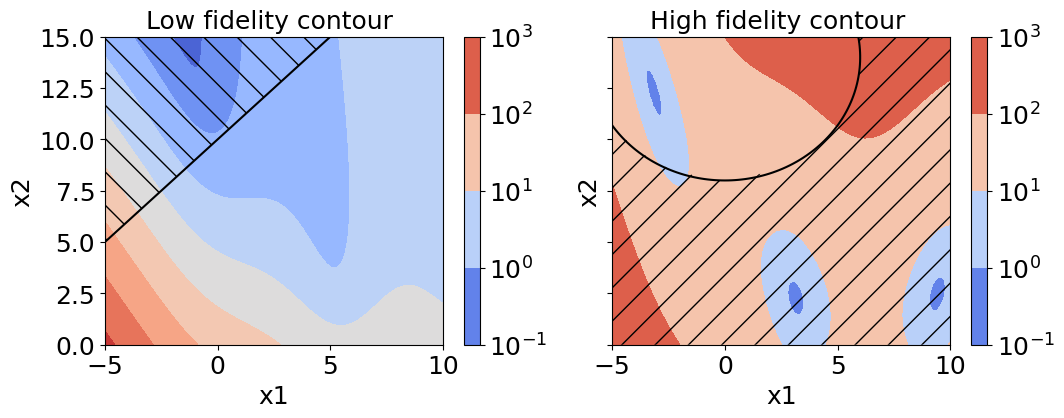}
\caption{Low-  and high-fidelity contour plots of example 3 and 4.}
  \label{fig:ex-Branin-mod}
    \end{subfigure}
  \caption{Contour plots of example 1 to 4. Feasible regions are indicated by the absence of parallel black lines.}
\end{figure*}

The objective of the first example is the Branin function~\citep{gardner2014,gramacy2016modeling,picheny2016bayesian,ariafar2019admmbo}, and its high-fidelity constraint function defines a circle. 
The low-fidelity constraint function is designed to resemble the high-fidelity one. 
The global minima of this example occurs at point $(-\pi,12.275)$ with minimum value $0.397887$. Furthermore, it is a feasible point with respect to the low fidelity problem. 
We deploy AECI as the high-fidelity acquisition function and  CUCB as the low-fidelity acquisition function.  
For the second example, we keep the same objective and constraint functions, while changing the low-fidelity acquisition function in the CMFBO algorithm to AECI. These two examples are designed to gain insights into the impact of changing low-fidelity acquisition functions by comparing the optimization paths of example 1 and 2.
Example 1 and 2 are visualized in Fig.~\ref{fig:ex-Branin}.

Example 3 and 4 also adopt the Branin objective function.
However, unlike example 1 and 2, the constraints in high- and low-fidelities are significantly different.   
Consequently, the feasible region of the low-fidelity problem does not contain the solution to the high-fidelity problem.
Thus, it is reasonable to expect reduced benefit from including additional low-fidelity samples. 
The optimal solution remains $0.397887$ at $(-\pi,12.275)$. Similar to example 1 and 2, we use AECI as the high-fidelity acquisition function and CUCB as the low-fidelity acquisition function for example 3 and AECI for both fidelities in example 4.  
Example 3 and 4 are visualized in Fig.~\ref{fig:ex-Branin-mod}.

The objective function of example 5 and 6 is the 2D Rosenbrock function~\citep{bananafunc}, also known as the banana function.
The high-fidelity constraint function $c_{h}$ defines the feasible region of the high-fidelity problem as a half-circle.
The constraints in high- and low-fidelities are similar and admit the same optimal solution.   
The global feasible minimum is $0$ at the point $(1,1)$ in both examples.
The optimal solution of the high-fidelity problem, \textit{i.e.},  $(1,1)$, is also the optimal solution to the low-fidelity problem. However, the feasible region of the low-fidelity problem is smaller, and its objective is scaled differently.
We use AECI as the high-fidelity acquisition function and CUCB as the low-fidelity acquisition function in example 5 and AECI for both fidelities in example 6.  
 
Example 7 and 8 employ the six-dimensional Hartmann function~\citep{park2017remarks} as the objective function. The feasible region of the low-fidelity constraint contains the optimal solution in the high fidelity. Given the increased dimensionality, we also increase the number of initial samples to 50 in example 8, to test their effect on the algorithm performance. 
 
We summarize the information for the synthetic problems in Table~\ref{tab:funcs}. For the constraint, we differentiate the correlation between fidelities by whether the optimal solution in high fidelity is feasible in the corresponding low fidelity. If so, the example is labeled `Y' in the table, otherwise `N'.

\begin{table*}[ht]
\centering
\begin{tabular}{llllll}
	\hline
	\hline
	Example &Objective&  $d$ &LF feasible & HF acquisition & LF acquisition  \\
	\hline
        1& Branin  &2 &Y &AECI  & CUCB\\
	\hline
        2& Branin  &2 &Y &AECI  & AECI\\
	\hline
        3& Branin  &2 &N &AECI  & CUCB\\
	\hline
        4& Branin  &2 &N &AECI  & AECI\\
	\hline
        5& Rosenbrock  &2 &Y &AECI  & CUCB\\
	\hline
        6& Rosenbrock  &2 &Y &AECI  & AECI\\
	\hline
    7& Hartmann  &6 &Y &AECI  & CUCB\\
	\hline
        8& Hartmann  &6 &Y &AECI  & AECI\\
	\hline
	\hline
\end{tabular}%
\caption{List of examples, the objective functions and dimension. LF feasible refers to whether low-fidelity feasible region contains the solution in high fidelity, where `Y' means yes and `N' means no. The acquisition functions used are listed as well. }\
\label{tab:funcs}
\end{table*}%

\begin{figure}
        \centering
	\includegraphics[width=0.98\linewidth]{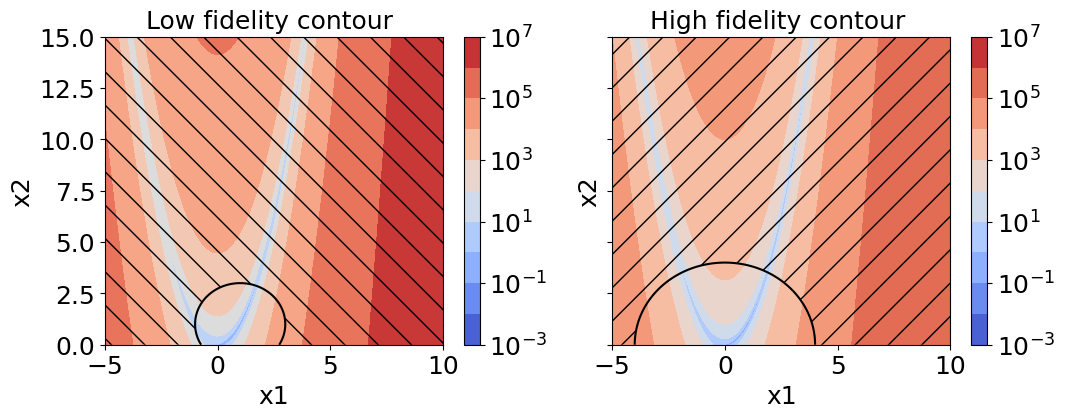}
  \caption{Low- (left) and high-fidelity (right) contour plots of example 5 and 6. Feasible regions are indicated by the absence of parallel black lines.}
  \label{fig:ex-banana}
\end{figure}

\begin{figure*}[ht]
    \begin{subfigure}{0.48\textwidth}
        \centering
	\includegraphics[width=\linewidth]{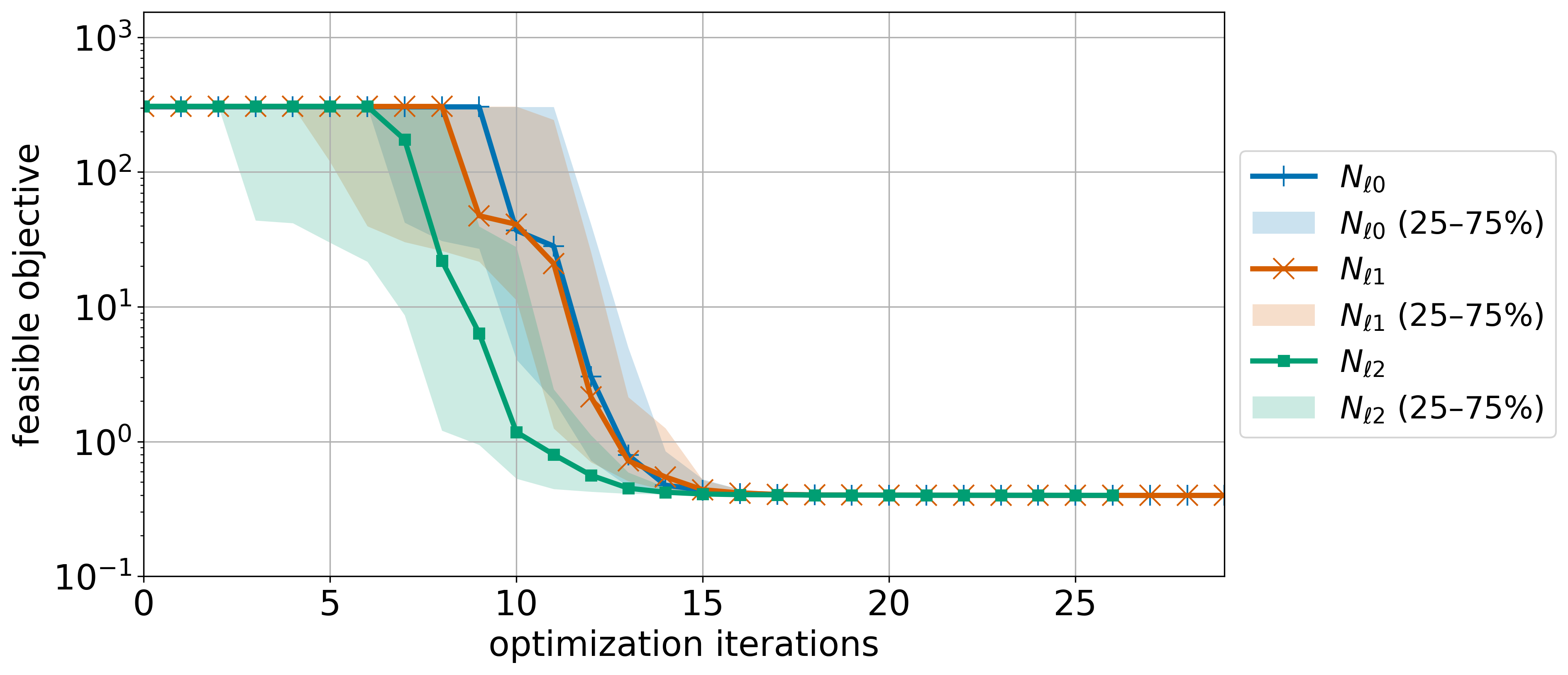}
  \caption{ Best feasible objective of example 1. }
  \label{fig:ex-Branin-sol}
    \end{subfigure}
    \hfill
    \begin{subfigure}{0.48\textwidth}
        \centering
	\includegraphics[width=\linewidth]{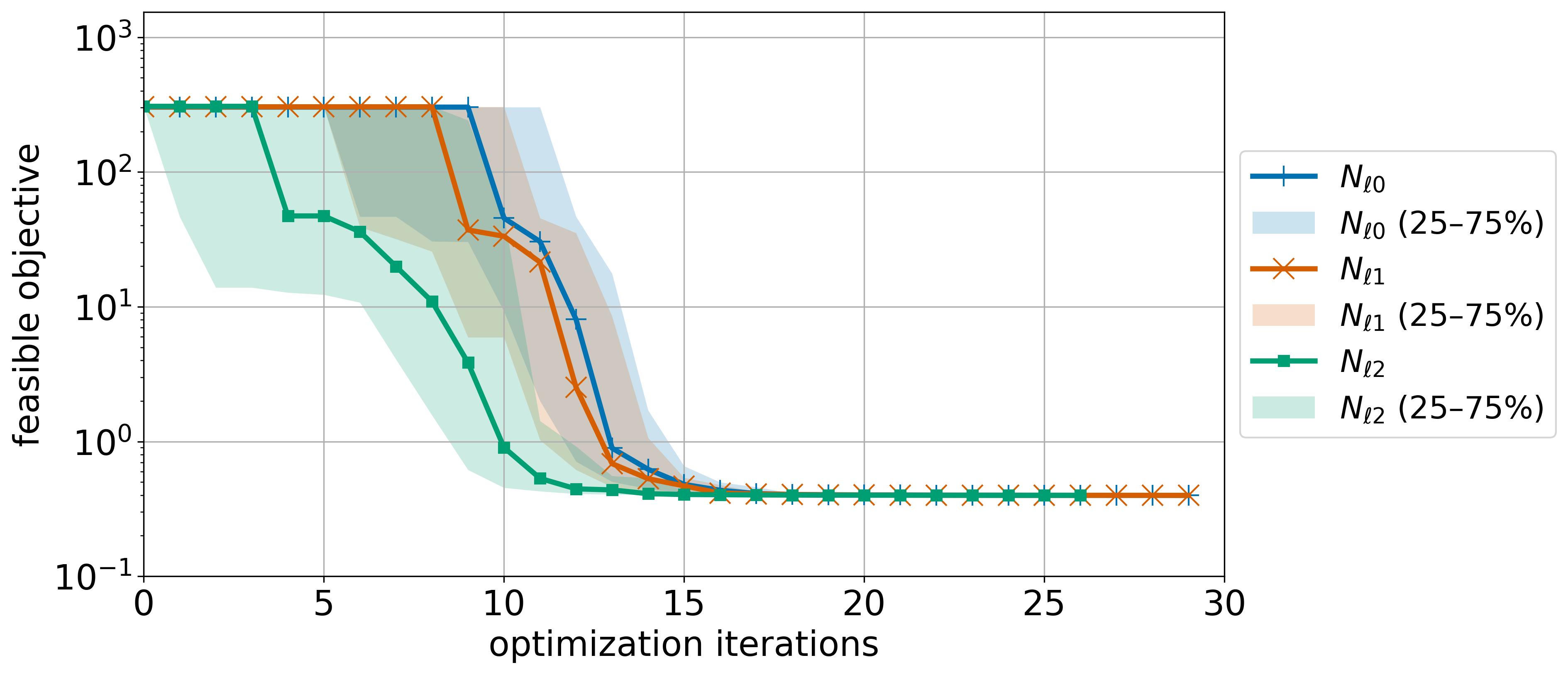}
  \caption{Best feasible objective of example 2.}
  \label{fig:ex-Branin-sol-aeci}
    \end{subfigure}
  \caption{CMFBO optimization history for example 1 and 2.}
\end{figure*}

For each example, the median best feasible objectives are plotted \textit{v.s.} the optimization iterations, along with shaded region representing the 25\% and 75\% of the best feasible objectives. 
In example 1 and 2, both the objective and constraint functions are highly correlated and the low-fidelity feasible region contains the optimal solution of the high-fidelity problem. Thus, it is more likely that low-fidelity samples can improve the accuracy of the multi-fidelity surrogate models. 
As illustrated in Fig.~\ref{fig:ex-Branin-sol} and~\ref{fig:ex-Branin-sol-aeci},
though all the experiments can quickly approach the correct optimum, a greater number of low-fidelity samples lead to more rapid convergence, taking advantage of the multi-fidelity data available.
Comparing example 1 and 2, using AECI as the low-fidelity objective improves the convergence performance of the CMFBO method. However, given that the solution is found quickly in both examples, the improvement is not significant.

\begin{figure*}[ht]
    \begin{subfigure}{0.48\textwidth}
        \centering
	\includegraphics[width=\linewidth]{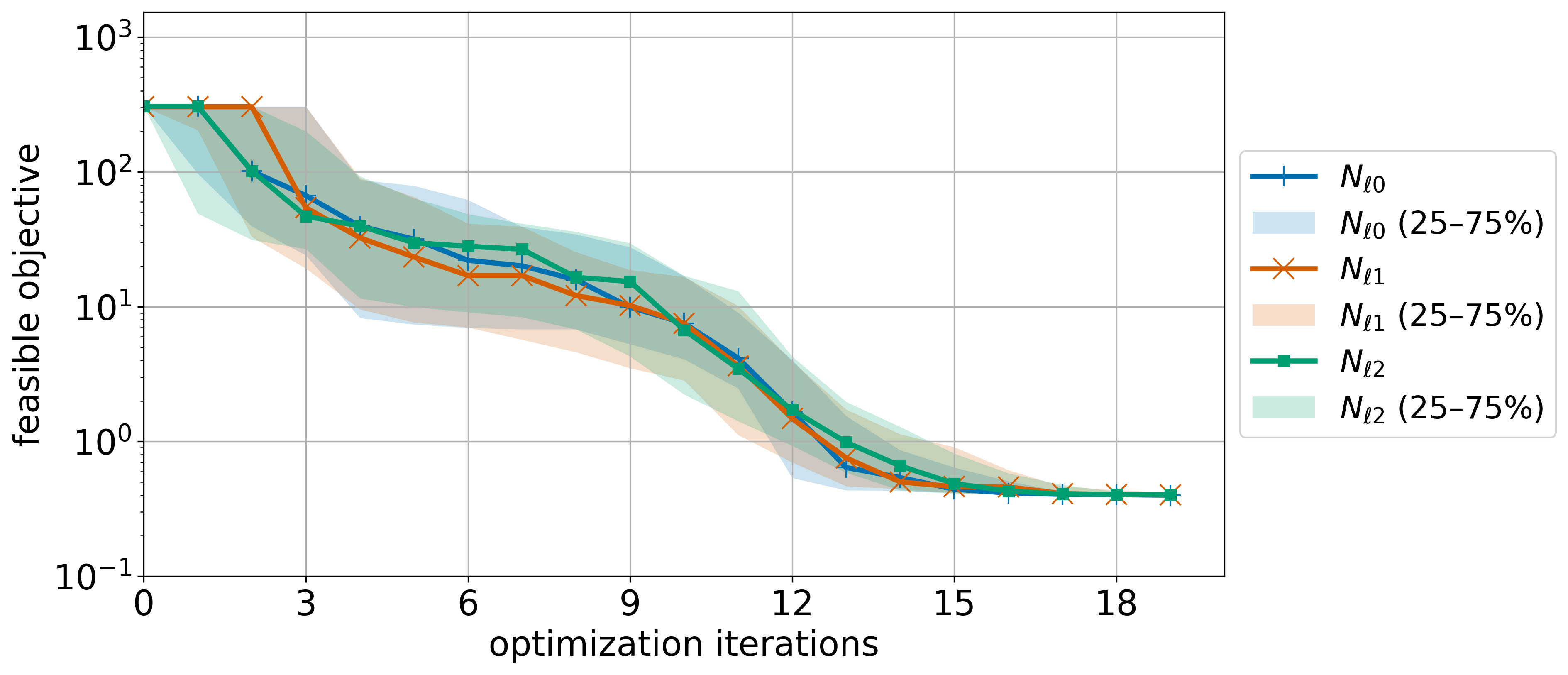}
  \caption{Best feasible objective of example 3. }
  \label{fig:ex-Branin-mod-sol}
    \end{subfigure}
    \hfill
    \begin{subfigure}{0.48\textwidth}
        \centering
	\includegraphics[width=\linewidth]{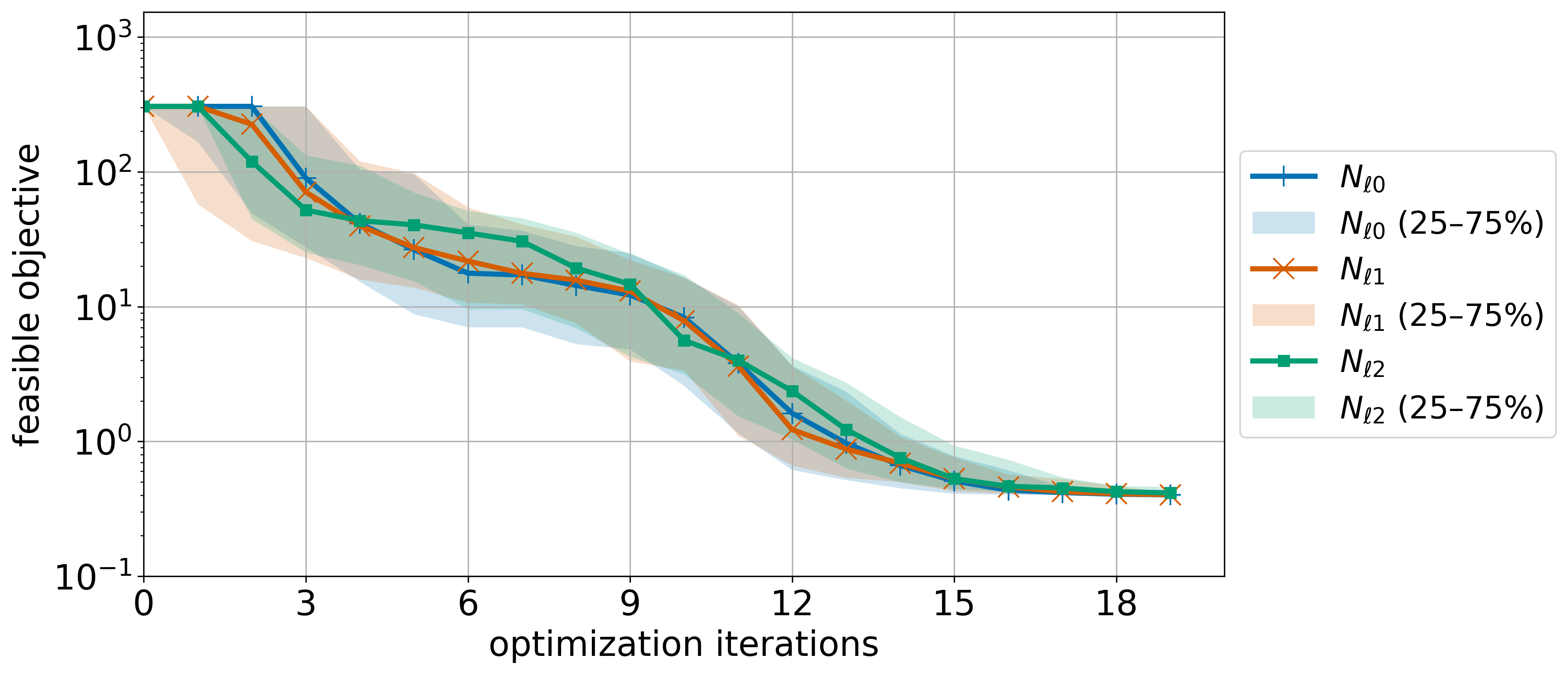}
  \caption{Best feasible objective of example 4.}
  \label{fig:ex-Branin-mod-sol-aeci}
    \end{subfigure}
  \caption{CMFBO solution for example 3 and 4.}
\end{figure*}

For example 3 and 4, the convergence history in Fig. \ref{fig:ex-Branin-mod-sol} and \ref{fig:ex-Branin-mod-sol-aeci} show that the initial feasible solution is identified quickly, within just two iterations, thanks to relatively large feasible regions. However, increasing the number of low-fidelity samples provides limited benefit in these two examples, due to their poor representation of the high-fidelity problem.
Comparing example 3 and 4, the difference due to the low-fidelity acquisition functions appears to be insignificant.

\begin{figure*}[ht]
    \begin{subfigure}{0.48\textwidth}
        \centering
	\includegraphics[width=\linewidth]{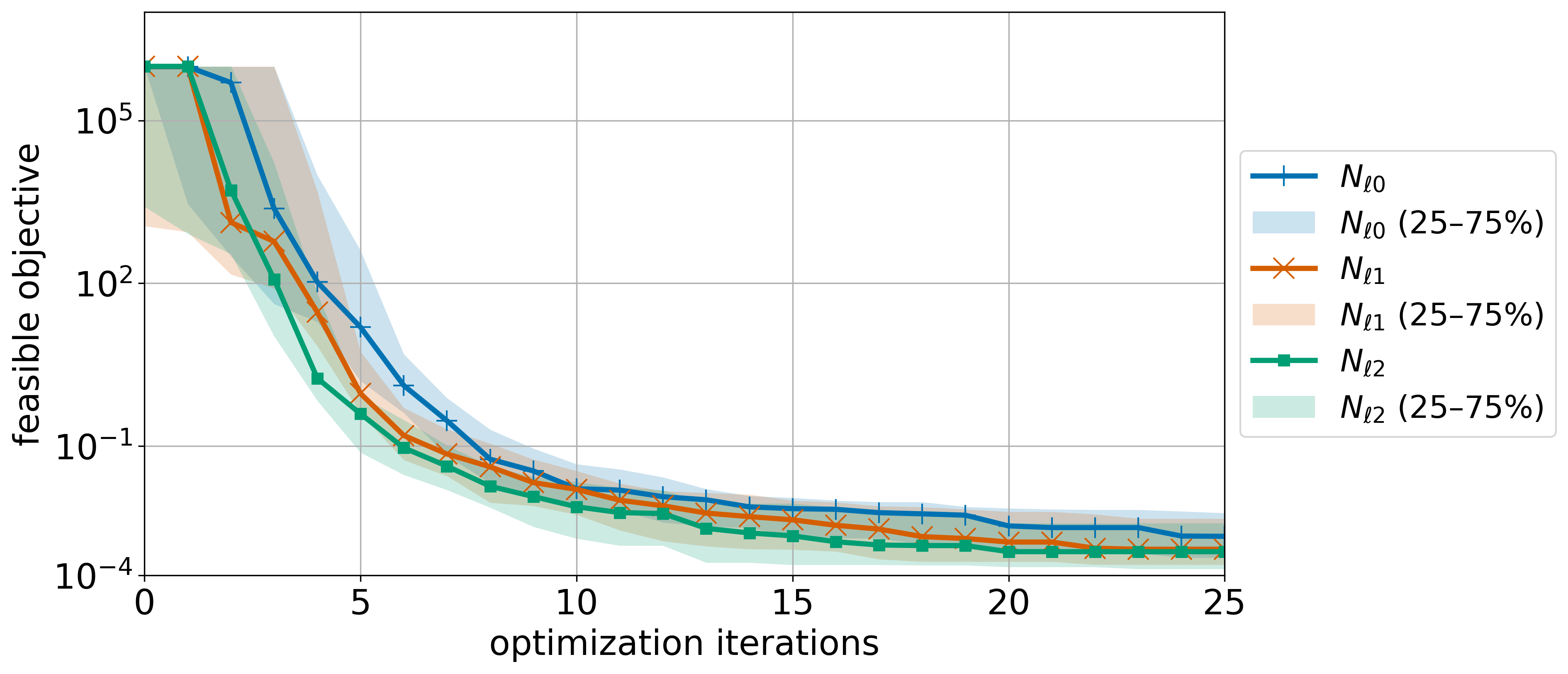}
  \caption{Best feasible objective of example
  5.}
  \label{fig:ex-Banana-sol}
    \end{subfigure}
    \hfill
    \begin{subfigure}{0.48\textwidth}
        \centering
	\includegraphics[width=\linewidth]{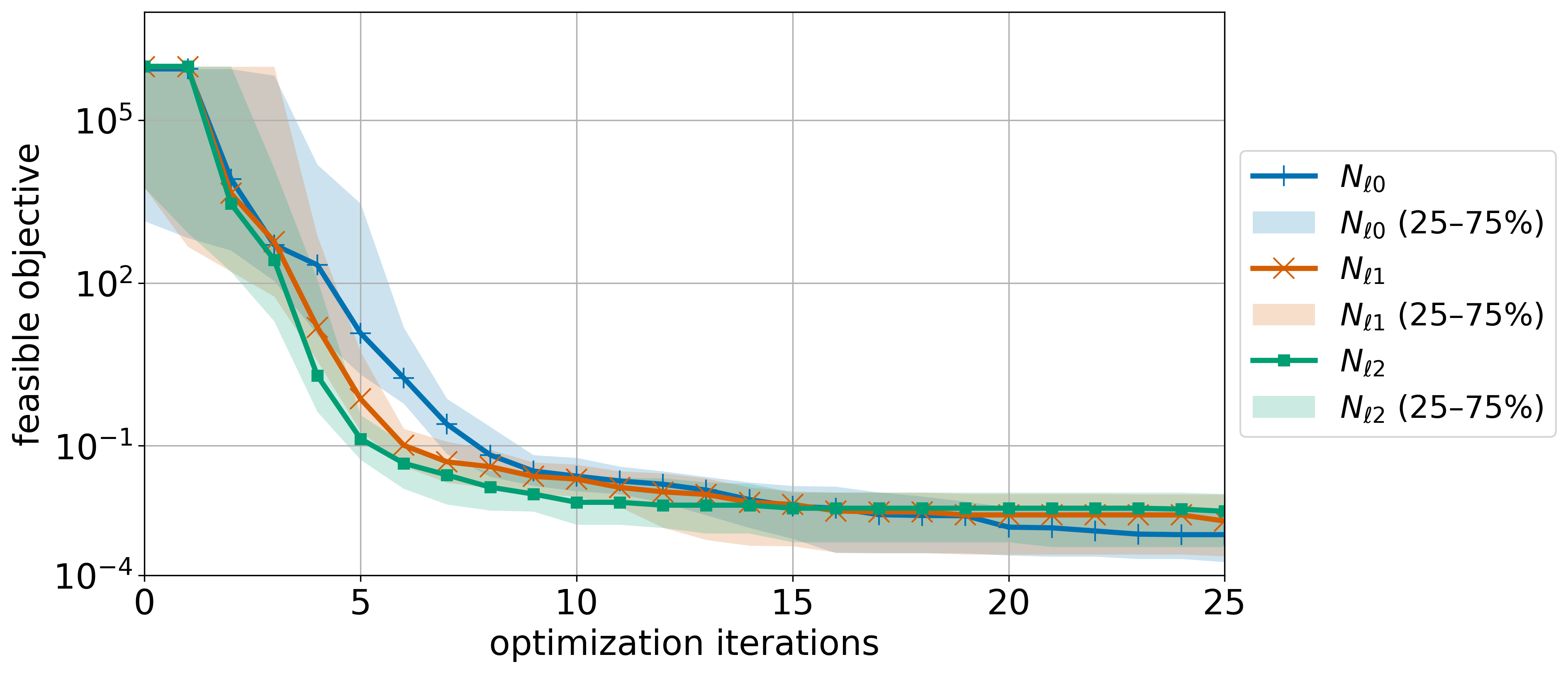}
  \caption{Best feasible objective of example 6.}
  \label{fig:ex-Banana-sol-aeci}
    \end{subfigure}
  \caption{CMFBO optimization history for example 5 and 6.}
\end{figure*}

 The best feasible objectives over the CMFBO iterations of example 5 and 6 are shown in Fig.\ref{fig:ex-Banana-sol} and~\ref{fig:ex-Banana-sol-aeci}, respectively. Fig.~\ref{fig:ex-Banana-sol} shows that all the experiments can quickly reduce the value of the best feasible objectives from $10^6$ to $10^{-3}$ within roughly 25 CMFBO iterations.
We highlight that the $N_{\ell 2}$ experiment takes 15 iterations to reach the similar accuracy obtained by $N_{\ell 0}$ with 24 iterations. This implies that increasing the number of low-fidelity samples per high-fidelity sample improves convergence performance of this problem.
Comparing example 3 and 4, CUCB offers an improvement in convergence as the low-fidelity acquisition function, particularly when $N_{\ell}$ is larger.
\begin{figure*}[ht]
    \begin{subfigure}{0.48\textwidth}
        \centering
	\includegraphics[width=\linewidth]{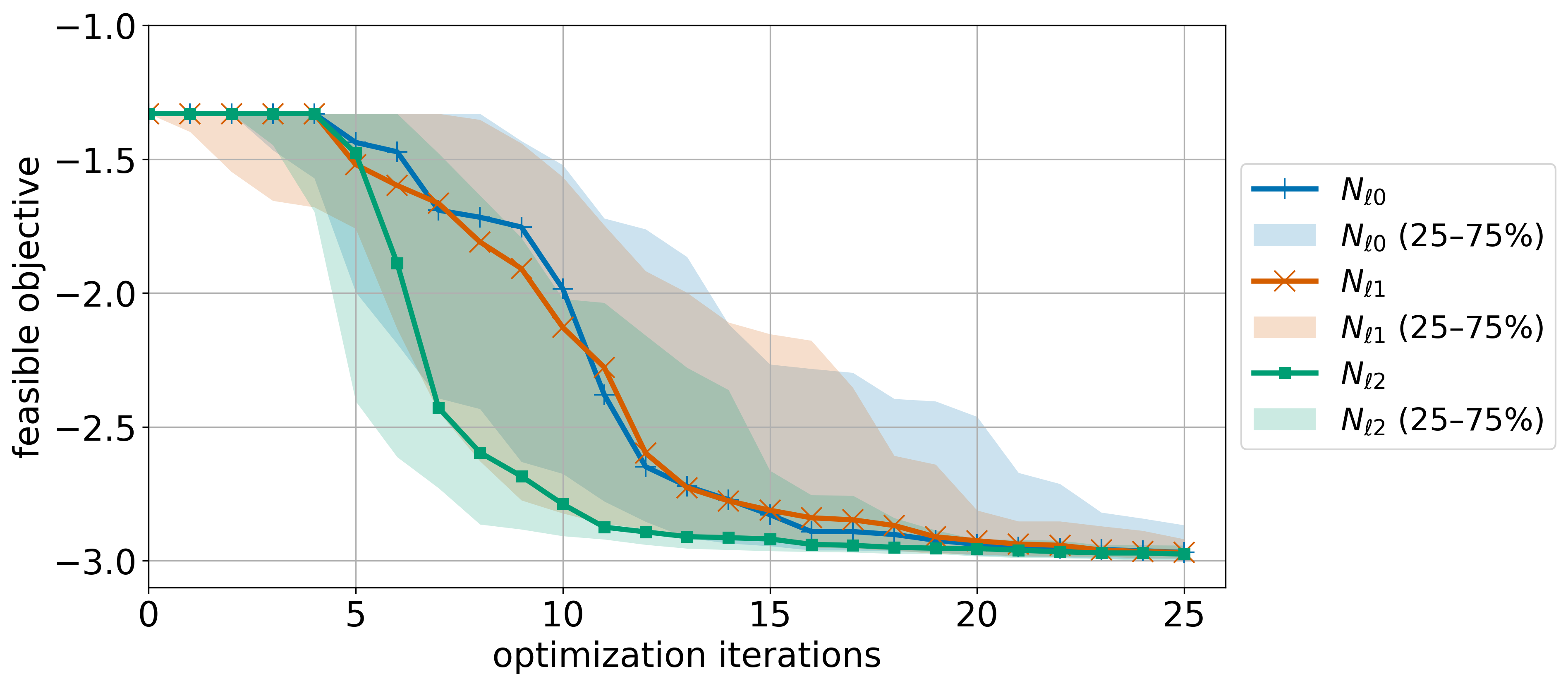}
  \caption{Best feasible objective of example  
  7.}
  \label{fig:ex-hartmann-aeci-ucb}
    \end{subfigure}
    \hfill
    \begin{subfigure}{0.48\textwidth}
        \centering
	\includegraphics[width=\linewidth]{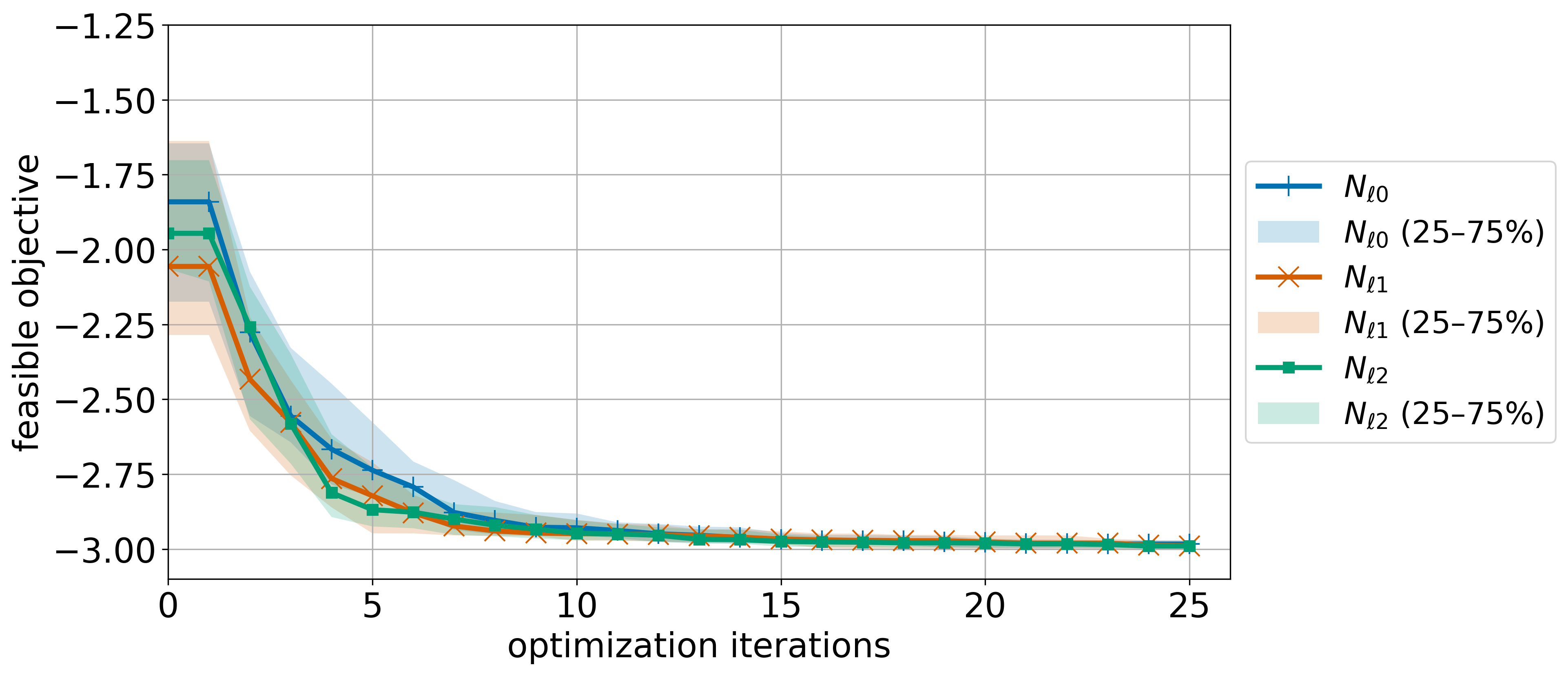}
  \caption{Best feasible objective of example 8.}
  \label{fig:ex-hartmann-aeci-aeci}
    \end{subfigure}
  \caption{CMFBO solution for example 7 and 8.}
\end{figure*}

For the six-dimensional Hartmann~\citep{park2017remarks} problem, Fig.~\ref{fig:ex-hartmann-aeci-ucb} and~\ref{fig:ex-hartmann-aeci-aeci} illustrate an effective best feasible objective reduction for all tests. Similar to previous examples, a larger $N_{\ell}$ slightly improves the convergence behavior. Increasing number of initial samples leads to a faster convergence and considerably reduces the uncertainty of the optimization runs at the cost of more initial sampling of the black-box functions. 

Our results demonstrate that low-fidelity samples improve the performance of CMFBO methods as long as the low-fidelity data preserve important features of the high-fidelity counterpart. Further, the best choice of acquisition functions appears to be problem dependent, thus validating the development of more acquisition functions.  

\subsection{An ICF design optimization problem}\label{se:icf}
In this section, we apply the proposed CMFBO method to an inertial confinement fusion (ICF) design problem with two design variables.
The objective function is based on simulation data collected using HYDRA~\citep{metal2001}, a multi-physics simulation code developed at Lawrence Livermore National Laboratory.
For a practical study of the performance of the algorithm, the CMFBO algorithms do not directly conduct HYDRA simulations due to the complexity of the software and the computational cost of the simulations.
Rather, we use an archived HYDRA database of one-dimensional capsule simulations with design variations based on National Ignition Facility (NIF) shot N210808, which is the first ICF experiment to exceed the Lawson fusion condition and produce more than 1 MJ of energy~\citep{N210808-lawson}. 

The high-fidelity database consists of standard HYDRA simulations, while the low-fidelity simulations are considered ``burn off'', in which the nuclear cross section has been artificially reduced by a factor of 1000, effectively turning off any yield amplification from alpha particle deposition. 
This ``burn-off'' model, while adopted as a simple simulation tool to create a low-fidelity physics model, has real-world applications. The best design of a low-fidelity model can be used to speed up an experimental campaign.
From the multi-fidelity optimization point of view, it is believed that there exists some level of correlation between the response surfaces at different fidelities, but the optima may or may not align.
For more details, the readers are referred to the unconstrained multi-fidelity Bayesian optimization work~\citep{wang2023multifidelity}.

We select design variables `t\_2nd' and `sc\_peak', which are the timing of the second shock and the strength of the peak radiation drive (see~\cite{wang2023multifidelity} for a explanatory figure of the two variables). The design objective is to maximize the nuclear yield, or equivalently, to minimize its negative value. 
The two constraint functions are chosen to be $c_1(\cdot)$ = `adiabat' and $c_2(\cdot)$= `vImp'.
The implosion velocity `vImp' (km/sec) is the maximum inward velocity of the fusion fuel achieved during the implosion, and `adiabat' is a unit-less measure $>1$ of the entropy of the fuel at the time of peak implosion velocity. 
In one-dimensional capsule simulation, lower values of adiabat and higher values of velocity tend to achieve higher yield; however, more complex models show that additional performance-limiting physics can turn on if either the adiabat is too small or the velocity is too large. As such, it is often desirable to solve design optimization problems which constrain these functions.
Consequently, we impose the constraints $c_1(x)\geq 4.25$ and $c_2(x)\leq 350$.

Similar to the synthetic problems, CMFBO is applied to the ICF example problem for 100 repeated runs in order to report performance results that are less sensitive to randomness. The same algorithmic parameters are used and AECI is chosen as the high-fidelity acquisition function and CUCB as the low-fidelity acquisition function. Further, we choose $N_{\ell}=1$ as the number of low-fidelity samples per iteration. 

\begin{figure}
        \centering
	\includegraphics[width=0.98\linewidth]{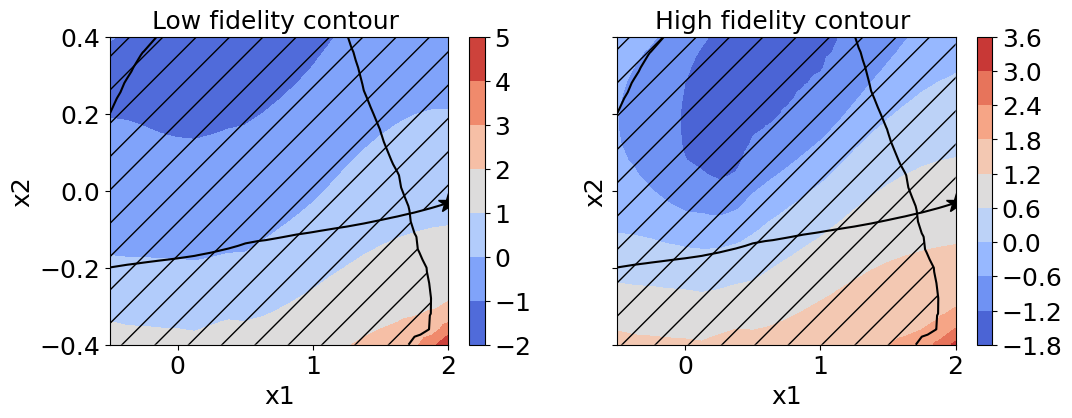}
  \caption{Low- (left) and high-fidelity (right) contour plots of the ICF design problem. The infeasible regions are marked with black dashed lines on the objective contour, along with the global optimum ($*$ sign).}
  \label{fig:ex-2DHydra}
\end{figure}
In Fig.~\ref{fig:ex-2DHydra}, we show contour plots (ground truth) of the scaled output nuclear yield over the 2D design space, obtained via exhaustive 100000 sample points. The objective plots in Fig.~\ref{fig:ex-2DHydra} clearly show that the unconstrained low- and high-fidelity simulation models predict different maximum outputs and distinct optimal design variables, highlighting the complexity of the multi-fidelity approach. However, the constraint functions are highly correlated. 

\begin{figure}
	\includegraphics[width=1\columnwidth]{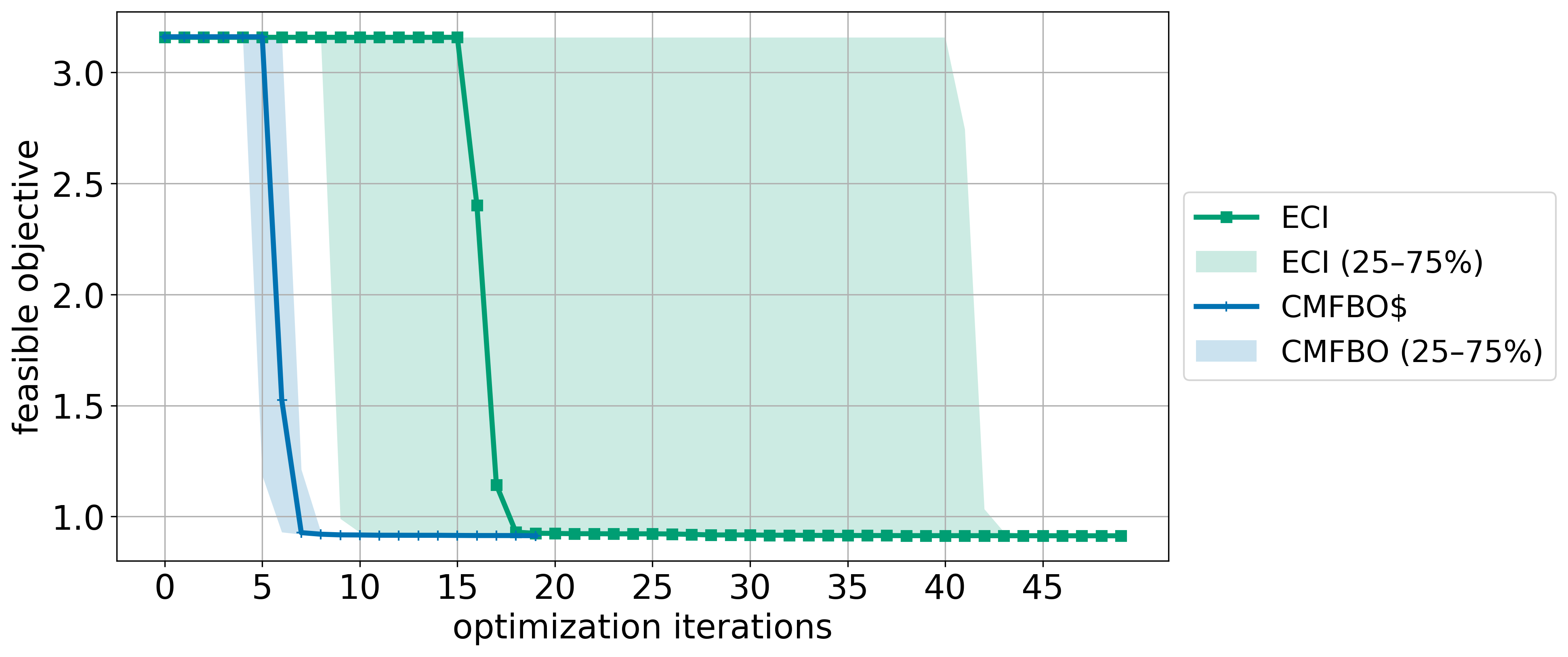}
  \caption{Example 2DHydra --- Best feasible objective from CMFBO with different parameters.}
  \label{fig:ex-2DHydra-sol}
\end{figure}

The best feasible objective curves from the proposed CMFBO are shown in Fig.~\ref{fig:ex-2DHydra-sol}.
In addition, we run the same problem with $N_\ell=1$, using ECI as the acquisition functions for both the high- and low-fidelity sampling.
When no feasible point exists, ECI uses random samples to continue. 
From Fig.~\ref{fig:ex-2DHydra-sol}, the proposed CMFBO method more quickly identifies an initial feasible point, converges faster, and exhibits better control over uncertainties thanks to the AECI acquisition function. Specifically, the 25th–75th percentile confidence interval is smaller, indicating improved stability of the proposed algorithm.


\subsection{A high-current joint design problem}\label{se:wedge}
\begin{figure}
  \centering
  \includegraphics[width=0.49\textwidth,trim={1cm 5cm 0cm 0cm},clip ]{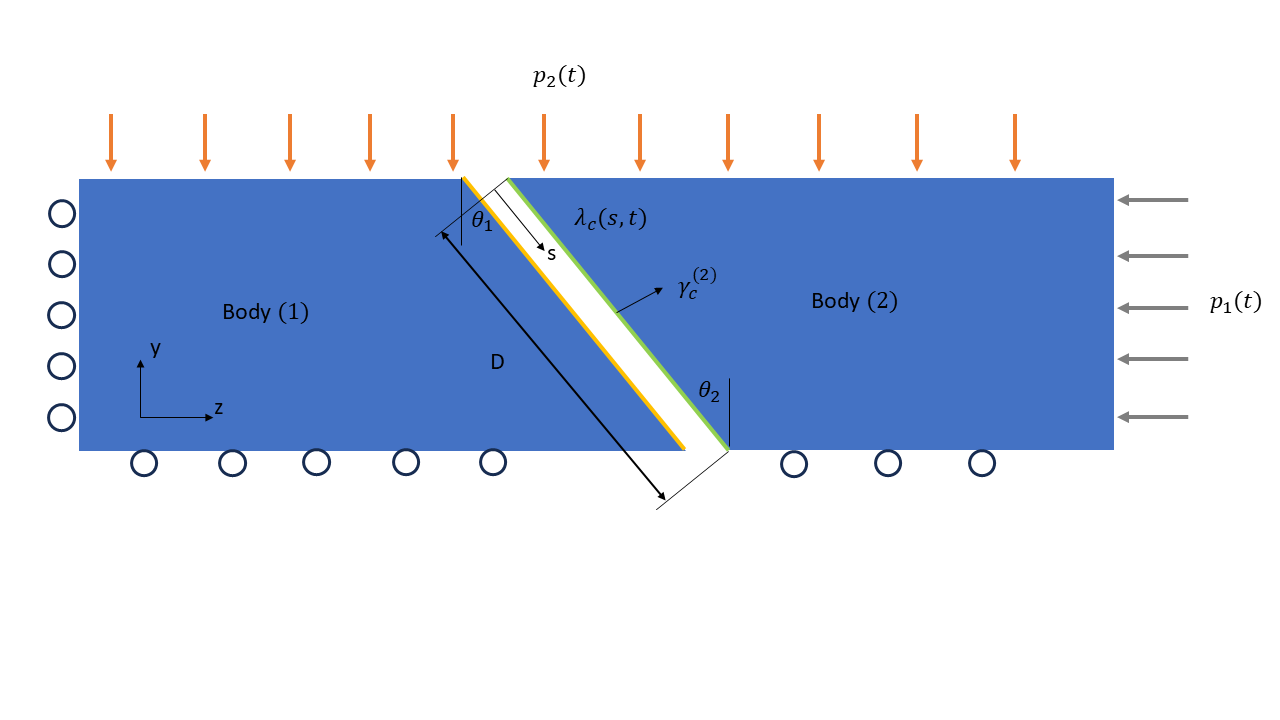}
	\caption{Wedge problem sketch in the $y-z$ plane. $p_1(t)$ and $p_2(t)$ are two time-dependent loads that stabilize at $P_1$ and $P_2$, respectively.}
\label{fig:wedge2d}
\end{figure}

This example is derived from the general problem of designing high-current joints, simulated using the finite-element method and computational contact mechanics.
A high-current joint is idealized as two wedges, subjected to external pressure loads $p_1(t)$ and $p_2(t)$. The geometry is shown in Fig.~\ref{fig:wedge2d}, where the current is assumed to flow roughly from the top surface and diffuse into the joint over time. The magnetic flux points out of the page, initially only in the ``air'' region above the top surface of the joint, but diffusing into the joint along with the current. Initially, the two halves of the joints are assembled with an initial preload force represented by the pressure $p_1(t)$, typically realized by clamping the two halves together and tightening an array of bolts (not shown).

Under the pressure load $p_1(t)$, the current is concentrated at the outer ($s=0$) region of the joint.   As the current builds up, a magnetic pressure is created, represented by the pressure load $p_2(t)$ exerted at the top of the wedges. 
Meanwhile, the current diffuses from the surface ($s=0$ at the outer fiber of the joint) to a depth $s=\alpha D$.
To prevent arcing or joint separation, it is important to maintain adequate contact pressure over the $s\in[0,\alpha D]$ surface, while not exceeding local material limits. 
Further, the maximum amount of preload $P_1$ should be maintained at a small value, as large values of contact preload require additional structure surrounding the joint. 
Therefore, the design optimization problem is formed to optimize the wedge geometry so that it maintains contact pressure within a working range as well as to minimize the maximum preload $P_1$.

\begin{figure}
  \centering
  \includegraphics[width=0.4\textwidth]{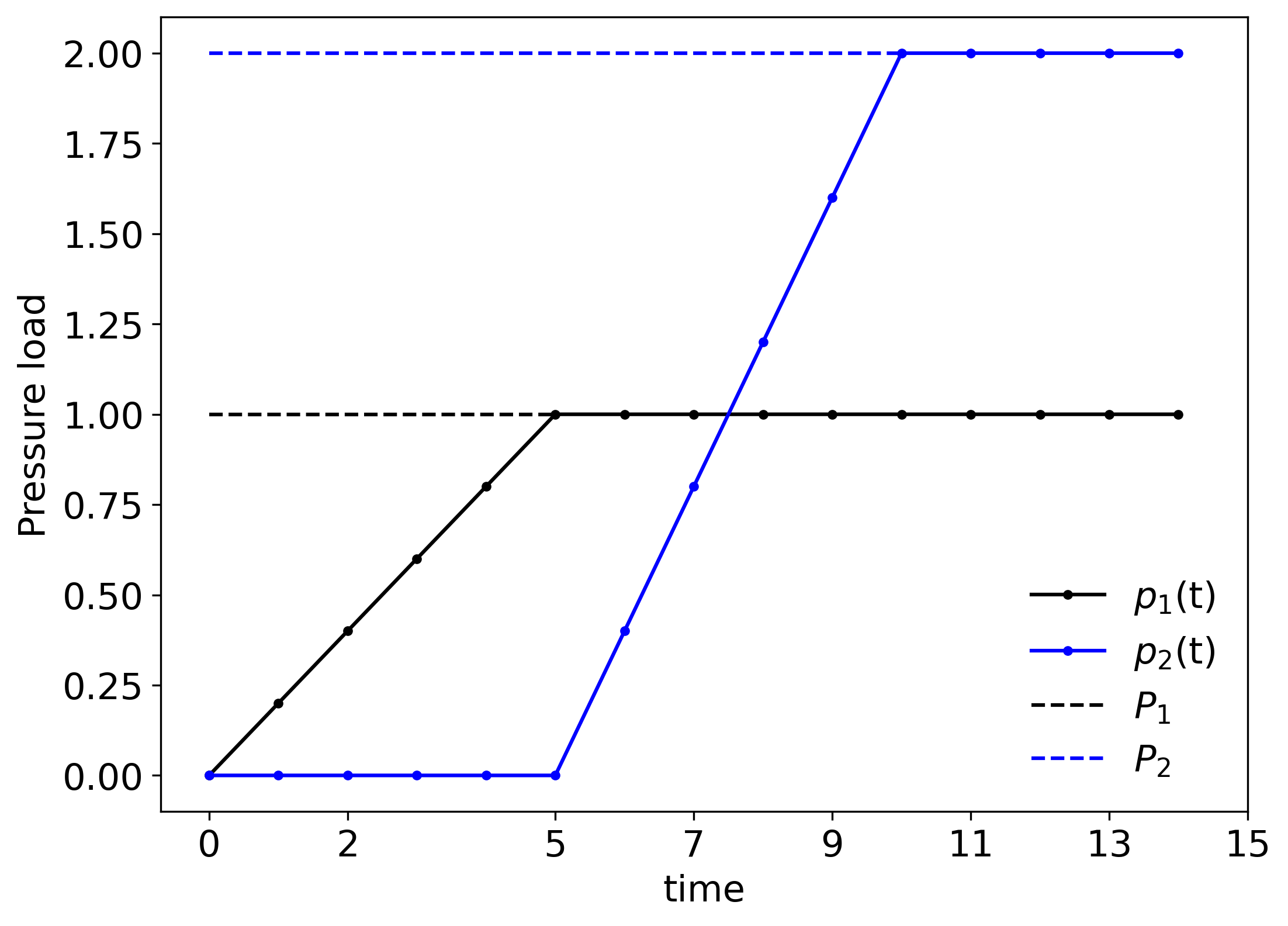}
	\caption{Illustration of the pressure load histories. The scalars $P_1$ (to be optimized) and $P_2$ are the maximum values of the loads.}
\label{fig:p1p2}
\end{figure}

The initial preload pressure $p_1(t)$ increases linearly from $0$ to a yet-to-be-optimized $P_1$ over the interval $[0,5]$, as seen in Fig.~\ref{fig:p1p2}; whereas the magnetic pressure $p_2(t)$ increases linearly from $0$ to $2$ over the interval $[5,10]$. The left face of the left wedge is fixed in the $z$ direction and the bottom faces of the both wedges are fixed in the $y-$direction, as illustrated in Fig.~\ref{fig:wedge2d}.

\begin{figure}
  \centering
  \includegraphics[width=0.5\textwidth,trim={6cm 15cm 6cm 10cm},clip] {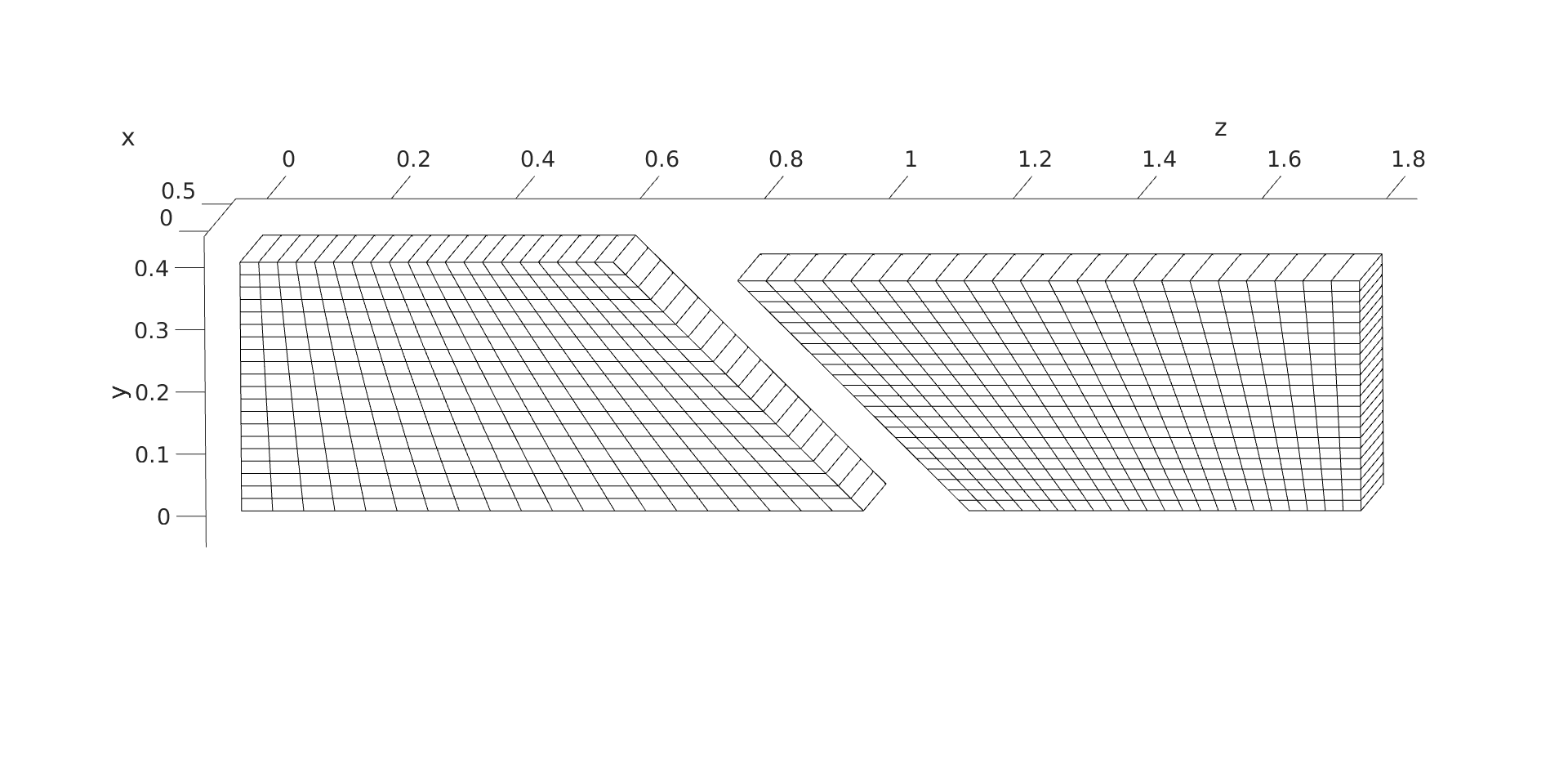}
	\caption{The undeformed initial design configuration of the wedges on the fine mesh.} 
\label{fig:wedge0}
\end{figure}

\begin{figure}
  \centering
  \includegraphics[width=0.5\textwidth,trim={6cm 15cm 6cm 10cm},clip] {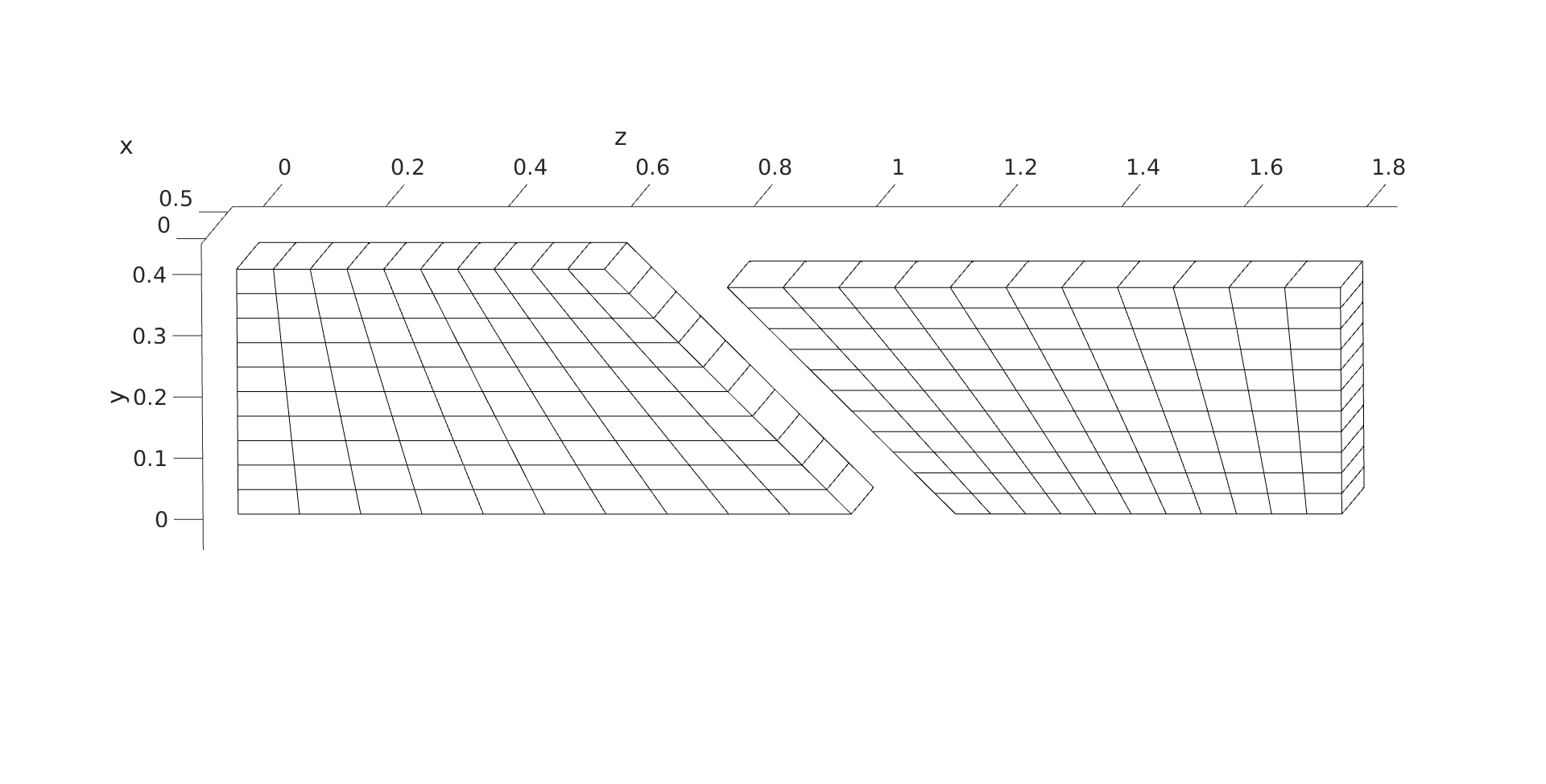}
	\caption{The undeformed initial design configuration of the wedges on the coarse mesh.} 
\label{fig:wedge-coarse}
\end{figure}

The wedges are discretized with linear hexahedron elements in two different mesh fidelities, as shown in Fig.~\ref{fig:wedge0} and Fig.~\ref{fig:wedge-coarse}. The refined mesh has quadruple the total number of elements of the coarse mesh and double the number of elements on the contact surface.
Furthermore, we impose plane strain boundary conditions on the front and back faces of the wedges by fixing displacement in the $x-$direction.
Due to the plane strain assumption, one layer of elements in the $x-$direction suffices.
Both wedges are modeled as an isotropic linear elastic material with the elastic modulus $E=200$ and Poisson's ratio $\nu=0.3$. 
The inclined surface on the right wedge is chosen to be the surface where the mortar contact integrals are performed, denoted as $\gamma_{c}^{(2)}$, where the superscript $(2)$ refers to the second body, and accordingly, $(1)$ refers to the wedge on the left.

The computational contact mechanics problem is solved by the finite element method and a mortar formulation~\citep{puso2004contact} wherein the gap function is integrated on the surface $\gamma_c^{(2)}$. 
The pressure constraints are imposed at time $5$ and $10$ to ensure the joint remains in desired working conditions given the time dependence of the two preloads.
Consequently, we solve the displacement and pressure fields at both time steps. 
We denote the nodal displacement field as $u_i^t, i\in I$, where $I$ is the collection of all node indices and $t\in \{5,10\}$. 
The nodal pressure field is denoted as $\lambda_i^t$, and is obtained via the solution of the contact mechanics problem. For more details of the finite element model and contact formulation, readers are referred to~\cite{wang2024design}.
We now present the mathematical formulation of the design optimization problem 
\begin{equation} \label{eqn:obj-1}
 \centering
  \begin{aligned}
   &\underset{\substack{x \in \Rbb^3}}{\text{minimize}} 
	  & & P_1 \\
   &\text{subject to}
	  && \lambda_{i}^{e,t} \geq p_l , \ s\leq \alpha D,\\
	   &&& p_u \geq \lambda_{i}^{e,t} , \ \forall s,\\
           &&& c(x)\geq 0.
  \end{aligned}
\end{equation}
The pressure constraint over the $s\in[0,\alpha D], \alpha=0.4$, region of the contact surface $\gamma_{c}^{(2)}$ is reflected by the lower bound $p_l=1$ on the contact pressure $\lambda_i^{e,t}$ (see below). 
In addition, as described above, an upper bound on the contact pressure $p_u=20$ is imposed over the contact surface. The additional constraints $c(x)$ are affine bound constraints on the optimization variable, which we specify later in~\eqref{eqn:ex1-bound}. 

\begin{figure}
  \centering
\begin{subfigure}{0.48\textwidth}
        \centering  
  \includegraphics[width=\linewidth,trim={0cm 6cm 0cm 4cm},clip]{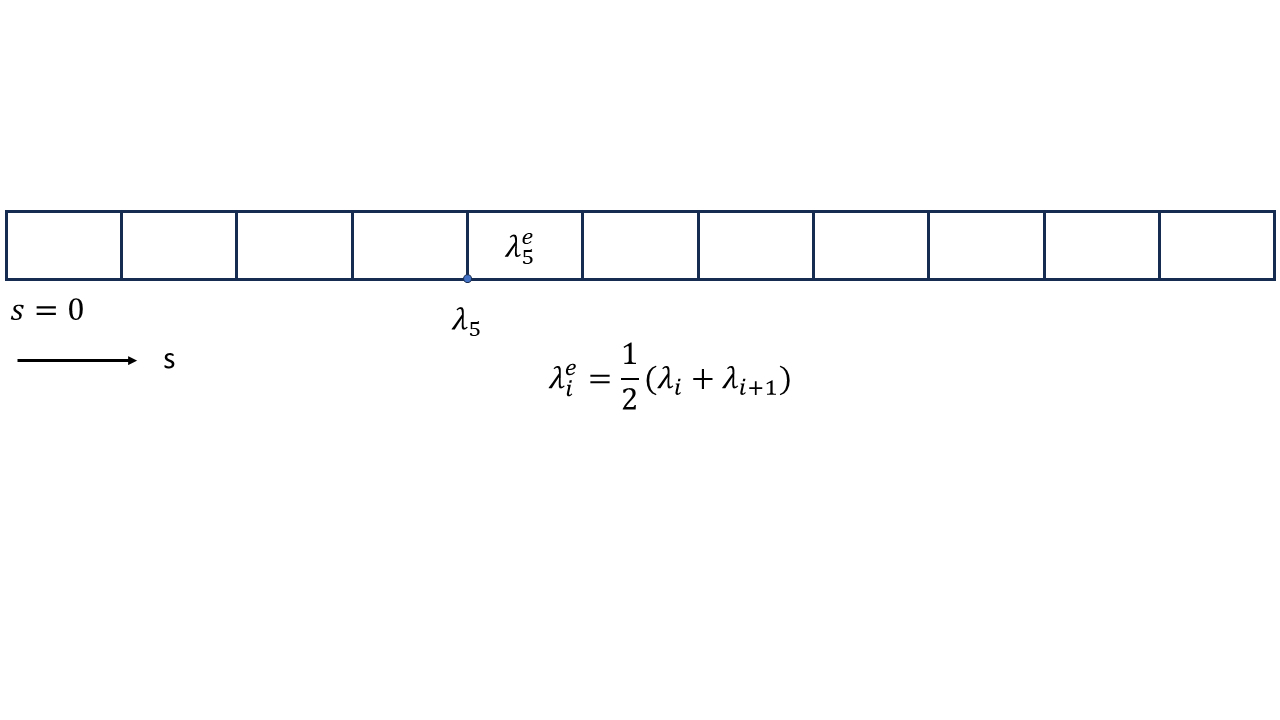}
	\caption{The element segment pressure example on the contact surface of the coarse mesh.}
\label{fig:elementseg_coarse}
\end{subfigure}
\begin{subfigure}{0.48\textwidth}
        \centering  
  \includegraphics[width=\linewidth,trim={0cm 6cm 0cm 4cm},clip]{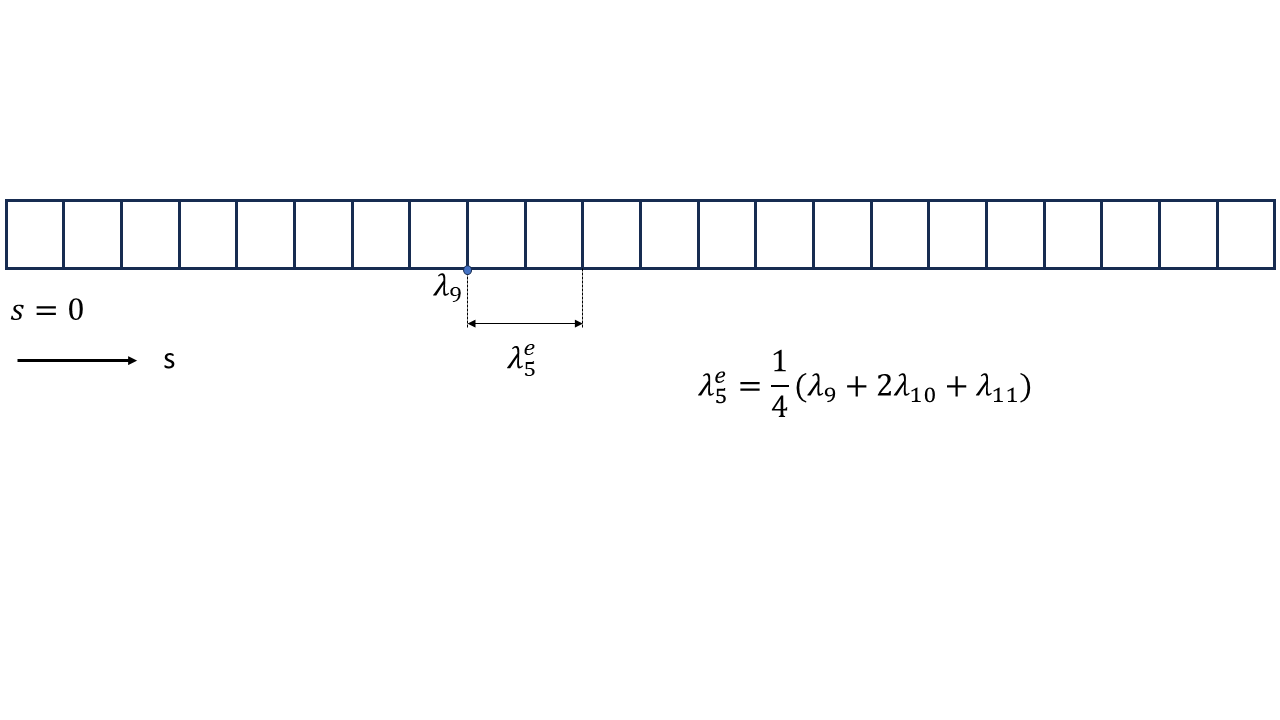}
	\caption{The element segment pressure example on the contact surface of the fine mesh.}
\label{fig:elementseg_fine}
\end{subfigure}
\caption{The element segment pressure example on the contact surface.}
\end{figure}

Design problems in contact can encounter significant nonsmoothness as the region in contact changes~\citep{hilding1999nonsmooth}. 
Further, the objective and constraints in both fidelities should provide comparable measures of the pressure of the same region. 
Thus, we impose the constraints in~\eqref{eqn:obj-1} on the \textit{element segment} pressure $\lambda_i^{e,t}$ (superscript $e$ for element) for the $i$th element segment, defined as the average nodal pressure on the element segment.
Each element segment contains one or more elements on the contact surface, as the name suggests. 
For the coarse mesh, each element is considered an element segment,  with a total of $11$ element segments on the contact surface, as illustrated in Fig.~\ref{fig:elementseg_coarse}. The first constraint in~\eqref{eqn:obj-1} is thus of dimension $4$ on the region $s\in[0,0.4 D]$ of $\gamma_{c}^{(2)}$.

For the refined mesh, we set one element segment to be equivalent to two elements. Since there are $22$ elements and $23$ nodes on the contact surface $\gamma_c^{(2)}$, we again have a total of $11$ element segments, as illustrated in Fig.~\ref{fig:elementseg_fine}. Therefore, the first constraint in~\eqref{eqn:obj-1} is imposed on a total of four element segments, equivalent to  eight elements, on $s\in[0,0.4 D]$  of $\gamma_{c}^{(2)}$.

The optimization variables $x$ are the wedge angles $\theta_1$, $\theta_2$ and the maximum preload pressure $P_1$. 
The additional affine and bound constraints in~\eqref{eqn:obj-1} are 
\begin{equation} \label{eqn:ex1-bound}
 \centering
  \begin{aligned}
	   &30^{\circ} \leq \theta_1,\theta_2\leq 60^{\circ},\\
	   &0.5\leq P_1 \leq 1.5,\\
	   &  \theta_1 \leq \theta_2 \leq \theta_1 + 8^{\circ}.
  \end{aligned}
\end{equation}
These constraints are imposed to ensure contact occurs at $t=5$ and $t=10$. For example, if the difference between the angles is too large, significant deformation can occur and cause the right wedge to slip out of contact with the left wedge.  

\begin{figure}
	\centering
         \includegraphics[width=0.5\textwidth,trim={6cm 15cm 6cm 10cm},clip]{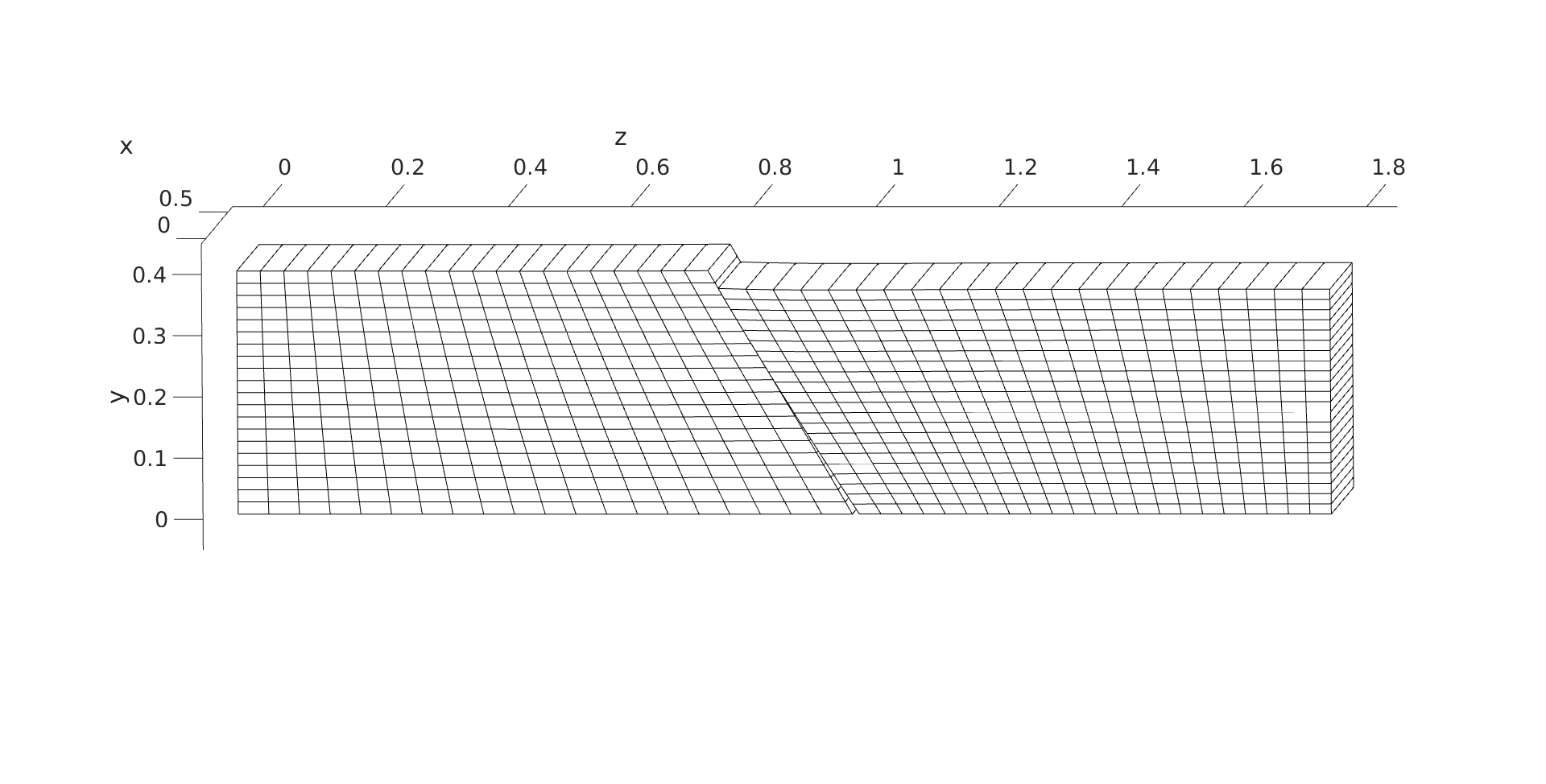}
	\caption{The deformed configurations of the optimal design at time 5$s$ (fine mesh). } \label{fig:wedge-5s}
\end{figure}

\begin{figure}
        \centering
         \includegraphics[width=0.5\textwidth,trim={6cm 15cm 6cm 10cm},clip]{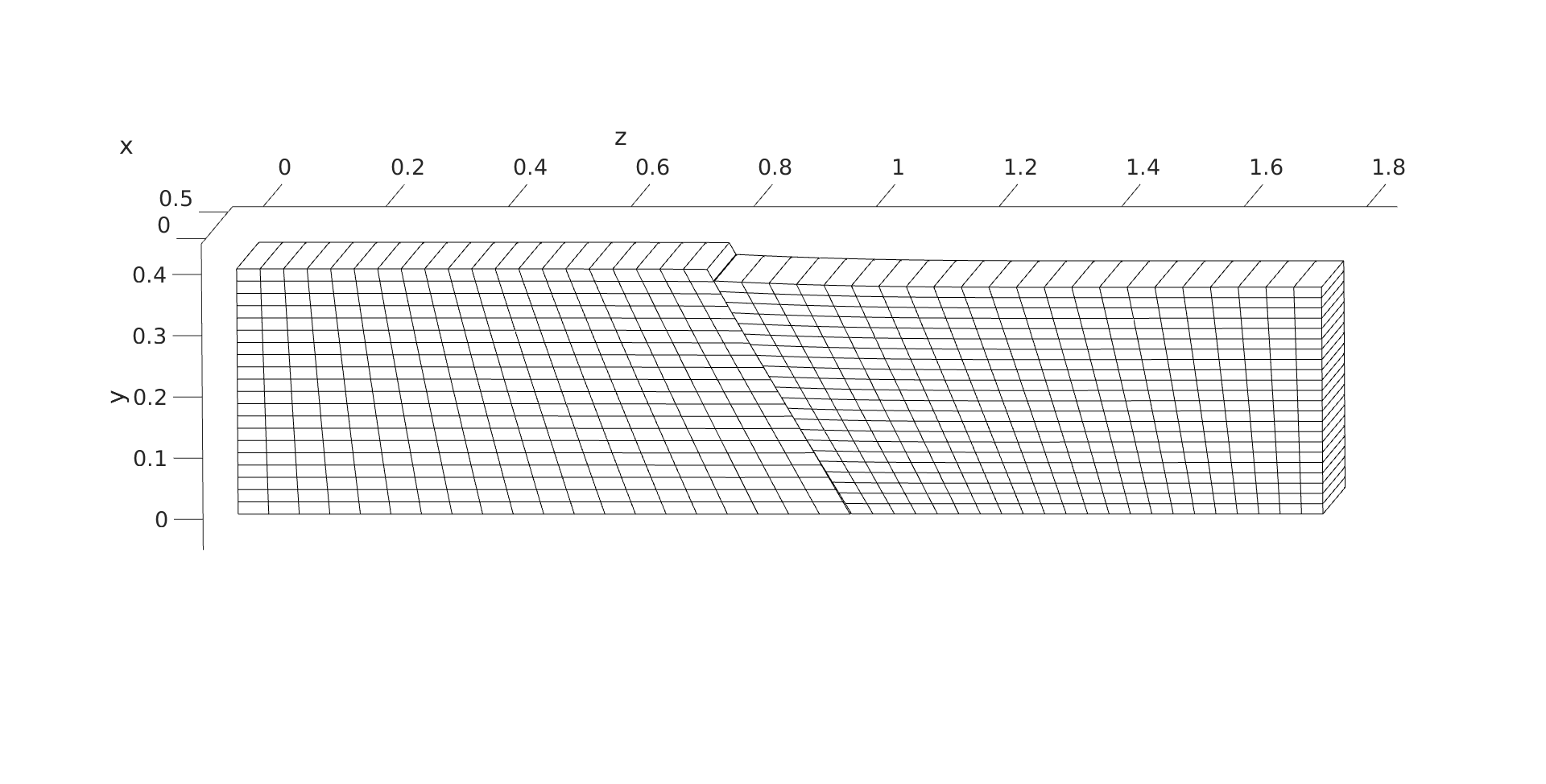}
	\caption{The deformed configurations of the optimal design at time 10$s$ (fine mesh). } \label{fig:wedge-10s}
\end{figure}

To better demonstrate the results of the CMFBO method, we also perform a gradient-based optimization with Ipopt~\citep{ipopt} for this problem, using the sensitivities of the objective and constraints with respect to the design variables. The sensitivities are computed by differentiating the fully discretized balance equations of the contact mechanics problem (see~\cite{wang2024design} for details).
The initial values of the variables are set to $\theta_1=45^{\circ}$, $\theta_2=45^{\circ}$ and $P_1=0.75$.
Using the gradient-based approach, the optimal solution is successfully found in 50 iterations and $57$ contact simulations  
with the values $\theta_1 =  30^{\circ}, \theta_2 = 31.9583^{\circ}$ and $ P_1 =  0.7527$.
The deformed configuration of the optimal design at times $5$ and $10$ is shown in Fig.~\ref{fig:wedge-5s} and~\ref{fig:wedge-10s}.

\begin{figure}
  \centering
  \includegraphics[width=0.48\textwidth]{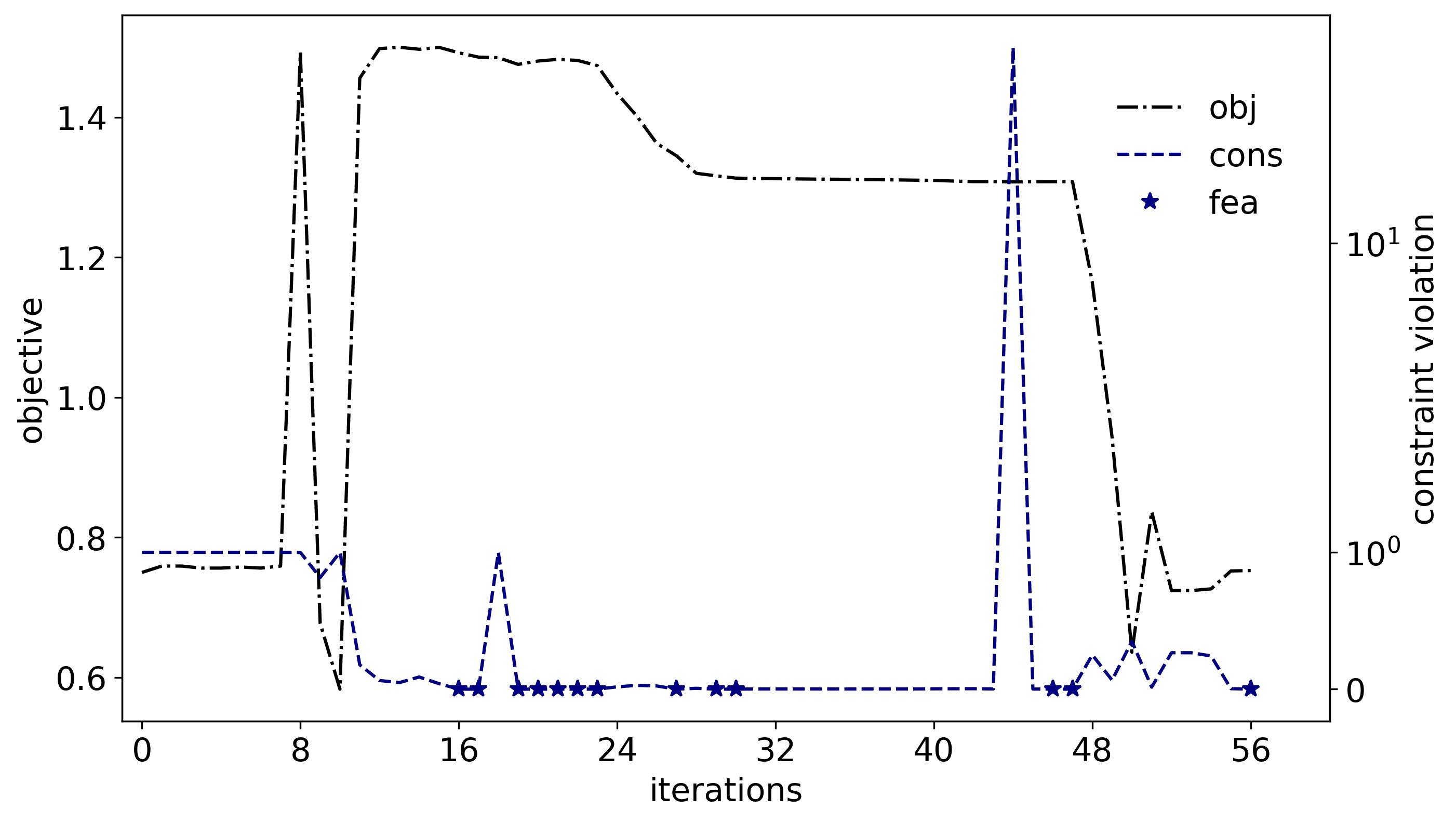}
	\caption{Gradient-based optimization design history for the wedge problem. The objective and constraint violation are plotted. The feasible solutions are marked by the $\star$ symbol on the constraint history curve. For better illustration, the constraint violation is plotted in symmetric log scale.}
\label{fig:wedge1}
\end{figure}

We show the optimization history in Fig.~\ref{fig:wedge1}. Due to the non-smoothness of the optimization problem, the objective versus iteration curve shows large oscillations and non-monotonicity. We mark the feasible iterates on the constraint curve and observe that only 13 of the 57 samples are feasible, due to the small feasible region of the problem. 

We apply the CMFBO method with AECI as both the high-fidelity and low-fidelity acquisition functions. The algorithmic parameters are set to $\alpha_0=1$, $N_f=5$, $N_\ell=1$ and $\beta=1$. Recall that a larger $N_f$ means the acquisition function switches from EMI to ECI when more feasible samples are found. 
The choice of a larger $N_f$ is due to the anticipated small feasible region.
Six initial designs are randomly chosen via Latin hypercube sampling in the three dimensional design space.
The CMFMO algorithm is run for $50$ iterations and repeated 5 times to obtain a median performance. The relatively low number of repeated runs is a result of the high computational cost of contact mechanics simulation. In addition, we apply CMFBO method under the same setup but with ECI as its acquisition functions in both fidelities. When no feasible initial samples are available, we again use random sampling to start ECI.

Despite the small feasible region of the problem, the proposed method succeeded in finding feasible samples at every run, with varying best feasible objectives. 
Using the Ipopt solution as the ground truth, Fig.~\ref{fig:wedge1} shows the error of the best feasible objective for both CMFBO, ECI, and gradient-based approaches.
We plot both the median and the best runs here given the small number of repeated runs and the fact the results of such a run would likely be adopted as the optimal solution to the design optimization problem. 
We cap the largest best feasible objective error as the the error between the maximum of the objective and the optimal value, \textit{i.e.}, $1.5-0.7527$.
The best CMFBO run obtains solution close to the ground truth design, and the median CMFBO run also finds feasible design with reduced objective in a comparable number of simulations. From Fig.~\ref{fig:wedge1}, the proposed method and acquisition functions continue to outperform ECI, whose median run does not make any progress with 56 samples. 
This is because even with 56 random samples, feasibility of the samples still cannot be guaranteed. While the gradient-based method used here can solve this problem, we emphasize that CMFBO does not require sensitivity derivation, computation, or implementation. 
\begin{figure}
  \centering
  \includegraphics[width=0.49\textwidth]{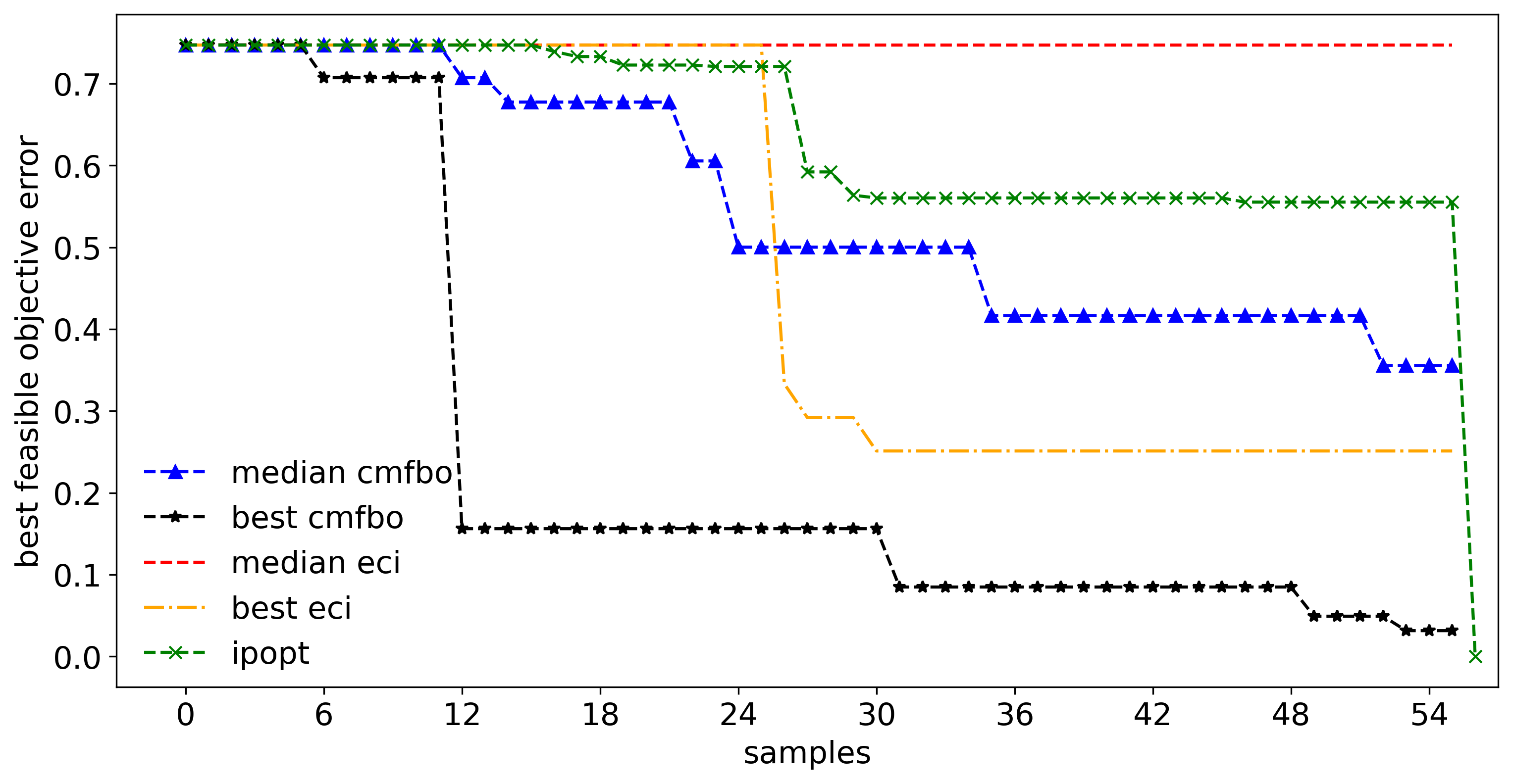}
	\caption{Best feasible objective error of CMFBO, ECI, and gradient-based optimization.}
\label{fig:wedge-cmfbo}
\end{figure}
%


\section{Conclusions} \label{sec:con}

We presented a novel constrained multi-fidelity Bayesian optimization method for design optimization problems with constraints. Our focus is on the design and implementation of novel acquisition functions in the CMFBO algorithm. With closed-form expressions, our method ensures ease of implementation and seamless integration into existing simulation software. Our proposed method overcomes the limitations of the widely used ECI acquisition function and provides promising alternatives. We demonstrate the effectiveness of our approaches through synthetic benchmark problems and two real-world applications.
Neither the constrained ICF nor the contact design optimization problem has previously been solved with Bayesian optimization. Therefore, our examples provide valuable validation of applying BO to identify feasible and desirable designs.

\appendix

\section{Numerical example}

For example 1 and 2, the high-fidelity problem has the following formulation 
\begin{align}  \label{eqn:ex-Branin-hf}
  \min_{x_1, x_2} \quad & f_h(x_1, x_2) := ( x_2 - \frac{5.1}{4\pi^2}x_1^2 + \frac{5}{\pi}x_1 - 6 )^2 \nonumber \\
  &\quad\quad\quad\quad\quad + 10\left(1 - \frac{1}{8\pi}\right)\cos(x_1) + 10 \\
  \text{s.t.} \quad & c_h(x_1, x_2) := 1.8-\nonumber\\
       &\sqrt{(x_1+2)^2+(x_2-12)^2} \geq 0, \\
              \quad & -5\leq x_1 \leq 10, \quad 0\leq x_2 \leq 15. 
\end{align}
  The global minima of example 1 and 2 in high fidelity occurs at point $(-\pi,12.275)$ and its minimum objective is $0.397887$.
The corresponding low-fidelity approximation of this problem is defined as
\begin{align}  \label{eqn:ex-Branin-lf}
  \min_{x_1, x_2} \quad & f_l(x_1, x_2) := 10\sqrt{f_h(x_1-2, x_2-2)} \nonumber \\
  &\quad\quad\quad\quad\quad + 2(x_1-2.5) - 3(3x_2-7) - 1 \\
  \text{s.t.} \quad & c_l(x_1,x_2):=1-\nonumber \\
           &\sqrt{(x_1+3)^2+(x_2-12.5)^2} \geq 0 , \\
              \quad & -5\leq x_1 \leq 10, \quad 0\leq x_2 \leq 15. 
\end{align}

The mathematical formulation of example 3 and 4 are given below.
\begin{align}  \label{eqn:ex-Branin-mod-hf}
  \min_{x_1, x_2} \quad & f_h(x_1, x_2) := \left( x_2 - \frac{5.1}{4\pi^2}x_1^2 + \frac{5}{\pi}x_1 - 6 \right)^2 \nonumber \\
  &\quad\quad\quad\quad\quad + 10\left(1 - \frac{1}{8\pi}\right)\cos(x_1) + 10 \\
  \text{s.t.} \quad & c_h(x_1, x_2) := 6- \\
              &\sqrt{(x_1)^2+(x_2-14)^2} \geq 0, \\
  \quad & -5\leq x_1 \leq 10, \quad 0\leq x_2 \leq 15. 
\end{align}
and 
\begin{align}  \label{eqn:ex-Banana-mod-lf}
  \min_{x_1, x_2} \quad & f_l(x_1, x_2) := 10\sqrt{f_h(x_1-2, x_2-2)} \nonumber \\
  &\quad\quad\quad\quad + 2(x_1-2.5) - 3(3x_2-7) - 1 \\
  \text{s.t.} \quad & c_l(x_1, x_2) := 10+x_1 - x_2 \geq 0, \\
  \quad & -5\leq x_1 \leq 10, \quad 0\leq x_2 \leq 15. 
\end{align}

The mathematical formulation of example 5 and 6 are given below.
\begin{align}  \label{eqn:ex-Banana-hf}
  \min_{x_1, x_2} \quad & f_h(x_1, x_2) := 100( x_2-x_1^2)^2 + (1-x_1)^2 \nonumber \\
  \text{s.t.} \quad & c_h(x_1, x_2) := 4-\sqrt{x_1^2+x_2^2} \geq 0, \\
  \quad & -5\leq x_1 \leq 10, \quad 0\leq x_2 \leq 15. 
\end{align}
The corresponding low-fidelity approximation of this problem is defined as
\begin{align}  \label{eqn:ex-Banana-lf}
  \min_{x_1, x_2} \quad & f_l(x_1, x_2) := 50(x_2-x_1^2)^2 + (1-x_1)^2\nonumber \\
  \text{s.t.} \quad & c_l(x_1, x_2) := 2-\nonumber\\
     &\sqrt{(x_1-1)^2+(x_2-1)^2} \geq 0, \\
  \quad & -5\leq x_1 \leq 10, \quad 0\leq x_2 \leq 15. 
\end{align}

The high-fidelity Hartmann 6 function for example 7 and 8  is 
\begin{equation} \label{eqn:hartmann-hf}
 \centering
  \begin{aligned}
		f(\mathbf{x})=&-\frac{1}{1.94}\left(2.58+\sum_{i=1}^{4} \alpha_{i} \right. \\
        & \quad\left. \exp \left(-\sum_{j=1}^{6} A_{i j} 
         \left(x_{j}-P_{i j}\right)^{2}\right) \right)\\
		\alpha=&(1.0,1.2,3.0,3.2)^{T}\\
		\mathbf{A}=&\left(\begin{array}{cccccc}
			10 & 3 & 17 & 3.50 & 1.7 & 8 \\
			0.05 & 10 & 17 & 0.1 & 8 & 14 \\
			3 & 3.5 & 1.7 & 10 & 17 & 8 \\
			17 & 8 & 0.05 & 10 & 0.1 & 14
		\end{array}\right)\\
		\mathbf{P}=&10^{-4}\left(\begin{array}{cccccc}
			1312 & 1696 & 5569 & 124 & 8283 & 5886 \\
			2329 & 4135 & 8307 & 3736 & 1004 & 9991 \\
			2348 & 1451 & 3522 & 2883 & 3047 & 6650 \\
			4047 & 8828 & 8732 & 5743 & 1091 & 381
		\end{array}\right)\\
		x_i \in& [0.1, 1], \ i = 1, \cdots, 6\\
		\mathbf{x}^* =& (0.20169,0.150011,0.476874,0.275332,\\
            & 0.311652,0.6573). 
  \end{aligned}
\end{equation}
The low-fidelity Hartmann 6 function is 
\begin{equation} \label{eqn:hartmann-lf}
 \centering
  \begin{aligned}
		&f(\mathbf{x})=-\frac{1}{1.94}\left(2.58+\sum_{i=1}^{4} \alpha_{i}' f_{exp}(\nu_i)\right) \text{and} \\
                &\nu_i= -\sum_{j=1}^{6} A_{i j}\left(x_{j}-P_{i j}\right)^{2} \\
		&\alpha'=(0.5,0.5,2.0,4.0)^{T}\\
		&f_{exp}(x) = \left( \exp\left(\frac{-4}{9}\right) + \exp\left(\frac{-4}{9}\right) \frac{x+4}{9}\right)^9.\\
  \end{aligned}
\end{equation}
The high-fidelity constraint function is designed to be a six dimensional ball. 
\begin{equation} \label{eqn:hartmann-con-hf}
 \centering
  \begin{aligned}
   c_h(x) := 0.25- \sum_{j=1}^6 (0.3-x_j)^2 \geq 0. \\
  \end{aligned}
\end{equation}
   The low-fidelity constraint function is a linear function  
\begin{equation} \label{eqn:hartmann-con-cf}
 \centering
  \begin{aligned}
   c_l(x) := 0.25 - \sum_{j=1}^6  \beta_j x_j \geq 0, \\
  \end{aligned}
\end{equation}
where $\beta= (0.1, 0.15, -0.17, 0.03, -0.01, -0.35)$.

\section*{Acknowledgments}
This work was performed under the auspices of the U.S. Department of Energy by Lawrence Livermore National Laboratory under contract DE--AC52--07NA27344 and the LLNL-LDRD Program under Project tracking No. 21-ERD-028.  Release number LLNL-JRNL-2004712-DRAFT. The data that support the findings of this study are available from the corresponding author upon reasonable request.

\section*{Replication of Results} 
The Python implementation of the proposed algorithm can be made available upon request by 
e-mailing the corresponding author. 
The algorithmic parameters used in the proposed method were performed as described
in the article.
\section*{Conflict of interest} On behalf of all authors, the corresponding author states that there is no conflict of interest.

\bibliographystyle{sn-basic}
\bibliography{bibliography}

\end{document}